\documentstyle[amssymb,amstex,graphicx,graphics,12pt,a4wide]{article}
\newtheorem{theorem}{Theorem}

\newtheorem{proposition}{Proposition}
\newtheorem{corollary}{Corollary}
\newtheorem{remark}{Remark}

\newcommand{\MaxEdges}{\operatorname{MaxEdges}}

\newcommand{\Vertices}{\operatorname{Vertices}}

\newcommand{\Area}{\operatorname{Area}}
\newcommand{\Cov}{\operatorname{Cov}}

\newcommand{\Var}{\operatorname{Var}}

\newcommand{\PLT}{\operatorname{PLT}}

\newcommand{\PVT}{\operatorname{PVT}}

\newcommand{\Segments}{\operatorname{Segments}}
\newcommand{\lparen}{\langle\langle}
\newcommand{\rparen}{\rangle\rangle}

\begin{document}

\title{Second-Order Properties and Central Limit Theory for the Vertex Process of Iteration Infinitely Divisible and Iteration Stable Random Tessellations in the Plane}
\author{Tomasz Schreiber\\
Faculty of Mathematics and Computer Science\\
Nicolaus Copernicus University, Toru\'n, Poland\\
\textit{e-mail: tomeks[at]mat.umk.pl}\\
Christoph Th\"ale\\
Department of Mathematics\\
University of Fribourg, Fribourg, Switzerland\\
\textit{e-mail: christoph.thaele[at]unifr.ch}}
\date{}
\maketitle

\begin{abstract}
The point process of vertices of an iteration infinitely divisible or more specifically of an iteration stable random tessellation in the Euclidean plane is considered. We explicitly determine its covariance measure and its pair-correlation function as well as the cross-covariance measure and the cross-correlation function of the vertex point process and the random length measure in the general non-stationary regime, and we give special formulas in the stationary and isotropic setting. Exact formulas are given for vertex count variances in compact and convex sampling windows and asymptotic
relations are derived. Our results are then compared with those for a Poisson line tessellation having the same length density parameter. Moreover, a functional central limit theorem for the joint process of suitably rescaled total edge count and edge length is established with the process $(\xi,t\xi)$, $t>0,$ arising in the limit, where $\xi$ is a centered Gaussian variable with explicitly known variance.
\end{abstract}
\begin{flushleft}\footnotesize
\end{flushleft}

\section{Introduction}

\begin{figure}[t]
\begin{center}
 \includegraphics[width=0.45\columnwidth]{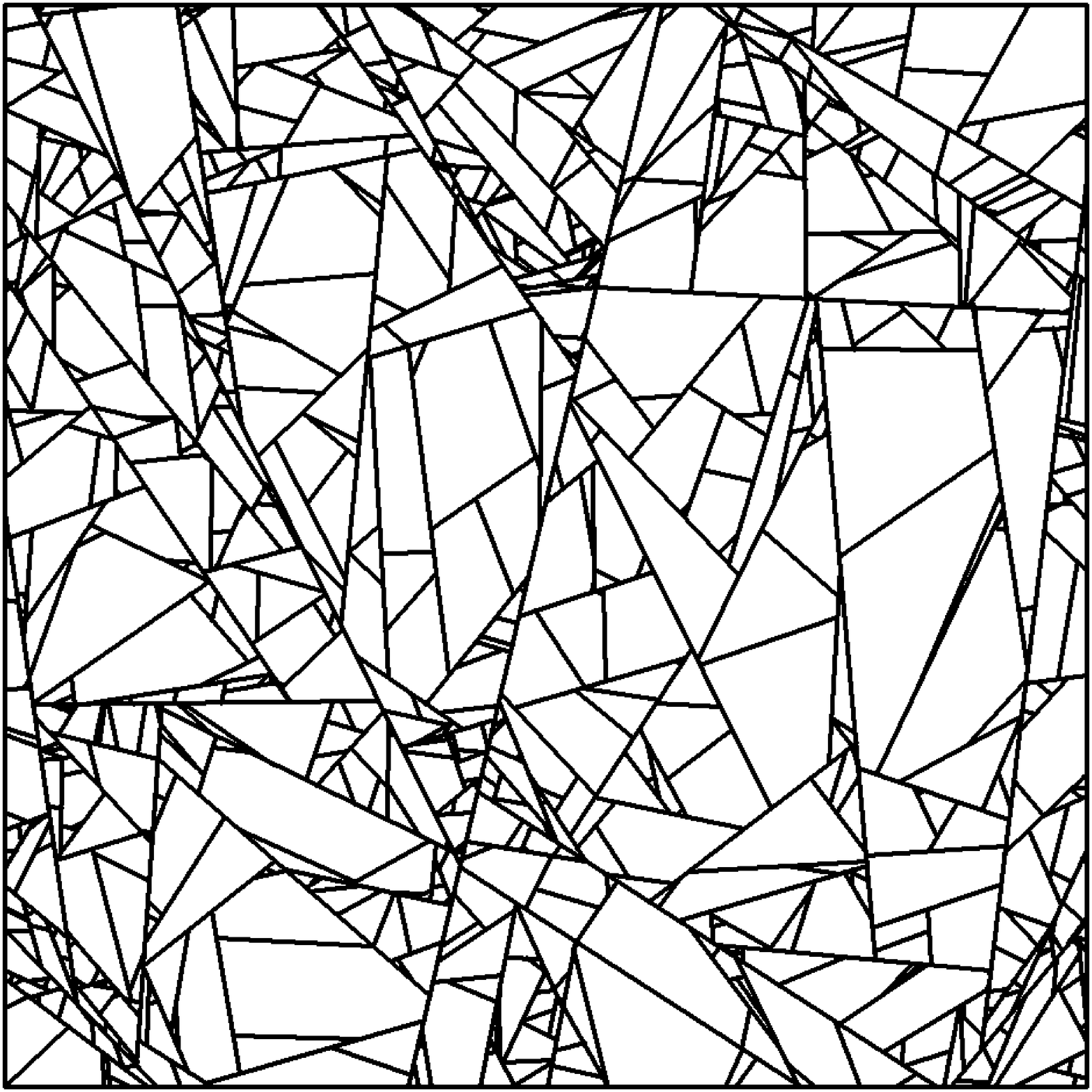}
 \includegraphics[width=0.45\columnwidth]{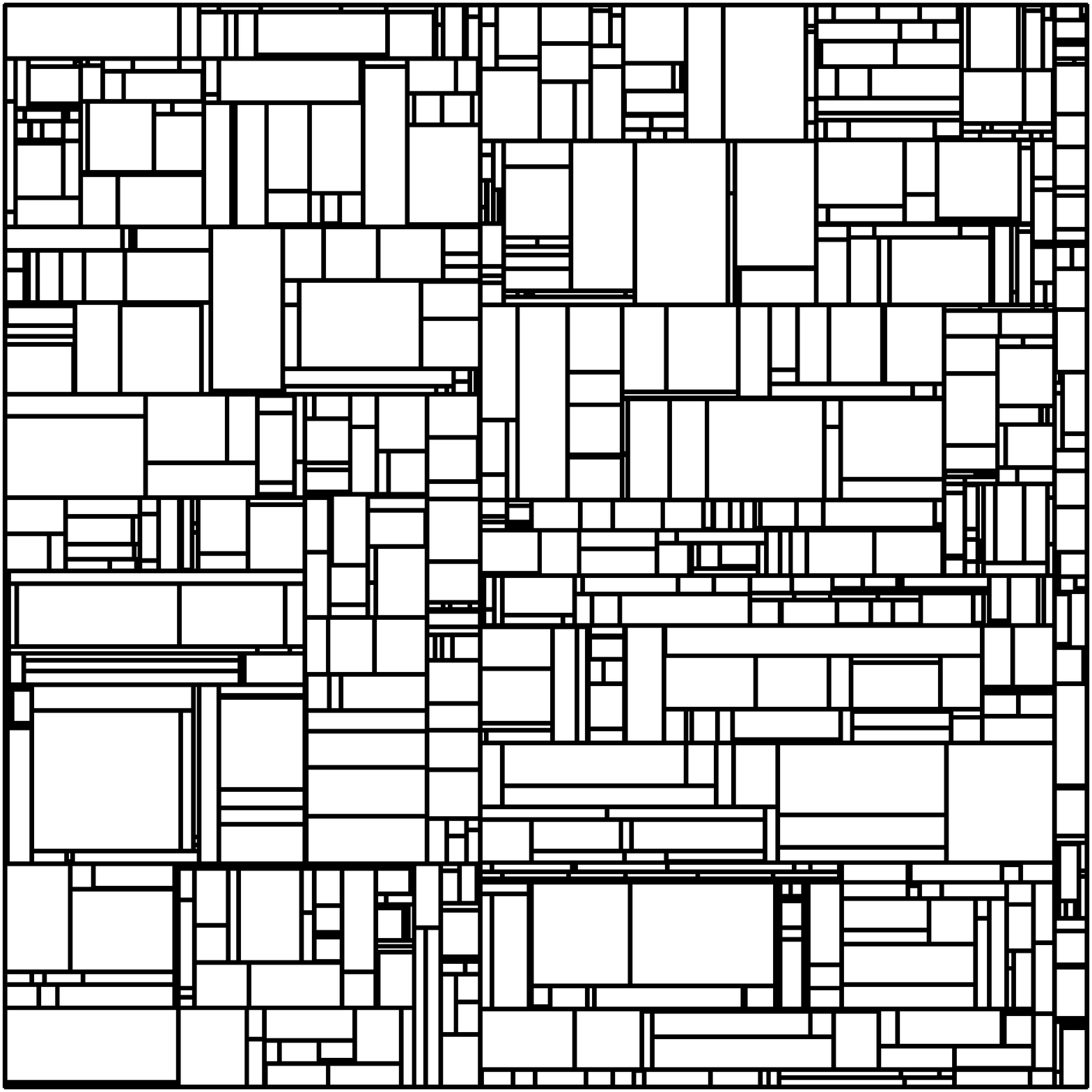}
 \caption{Realizations two STIT tessellations with different driving measures, the invariant measure (left) and a measure concentrated on lines pointing in only two orthogonal directions with weight $1/2$ (right)}\label{fig1}
\end{center}
\end{figure}

Random tessellations have attracted particular interest in stochastic geometry because of their wide applications ranging from classical geological problems to recent developments in telecommunication, see \cite{HS} and \cite{SKM}. It is one of the main purposes of the related theory to develop new classes of random tessellations that are mathematically tractable and yet rich enough in structure so that they may serve as new reference models for applications beside the classical Poisson hyperplane and the Poisson-Voronoi tessellation. A very recent model, the
so-called random STIT tessellations -- the abbreviation STIT stands for \textbf{st}able under \textbf{it}eration, see below -- was introduced by W. Nagel and V. Wei\ss\ in \cite{NW05}. One of the main features of these tessellations that distinguishes them from the above mentioned model classes is the property that their cells are not side-to-side, see Figure \ref{fig1}. This causes new geometric effects whose planar first-order properties in terms of mean values were explored in \cite{NW06}.\\
This paper will deal with second-order characteristics and central limit theory for the planar case, a topic that has been considered first by Wei\ss , Ohser and Nagel in \cite{NOW} and later also studied in \cite{ST}, where beside other characteristics the variance of the total {\it edge length} was determined. Here, in contrast, we will deal with second-order properties of the point process of {\it vertices} of the tessellation. To provide more general results, we will not restrict here to the class of stationary random iteration stable tessellations (STIT tessellations), and instead we will study the larger class of non-stationary random iteration infinitely divisible tessellations in the plane. They were introduced in \cite{ST} as generalizations of STIT tessellations and  in Section \ref{secstit} we will recall their construction in the spirit of Mecke-Nagel-Wei\ss\ (MNW) \cite{NW05}. It is an important observation that this spatio-temporal MNW-construction can be interpreted as a continuous time Markov process on the space of tessellations, whence the general theory of Markov processes is available. Before extending some mean value relations from the stationary iteration stable to the non-stationary iteration infinitely divisible case in Section \ref{secmv}, we will formulate the main technical tools from the theory of Markov processes in Section \ref{secmart}, on which our main results are based. They are the content of Sections \ref{secvar} and \ref{secclt}. The variance of the total number of vertices and that of the total number of maximal edges in a bounded observation window will be calculated in Section \ref{secvargen} for very general driving measures, whereas in Section \ref{secvariso} we specialize to the motion-invariant case. The vertex pair-correlation measure for general line measures is considered in Section \ref{secgenPCF} and the considerations are specialized again to the stationary and isotropic setting in Section \ref{secPCFiso}. Moreover, we determine in Section \ref{secgenEVPCF} the exact cross-covariance measure of the vertex point process and the length measure in the general case and provide specialized formulas in Section \ref{seccrosspcfiso} for the stationary and isotropic regime. Here, the most explicit formulas are available. Another topic treated there concerns the variance asymptotics for a sequence $W_R$ of growing windows as $R$ tends to infinity. Based on these results, Section \ref{secclt} deals with the central limit problem. It will be shown that a certain rescaled bivariate process of edge counts and edge lengths of a stationary random iteration stable (STIT) tessellation converges to the process $(\xi,t\xi)$, where $\xi$ is a centered normal random variable with an explicitly known variance.\\ We would like to emphasize that our results may be of interest in the context of statistical model fitting of random tessellations to real data, see \cite{HS}. Our central limit theory could be a base for statistical inference of tessellation models and related functionals, and asymptotic confidence intervals and statistical tests can be derived from them, since we make the first- and second-order moments explicitly available.\\ For general notions and notation related to stochastic geometry and the theory of random tessellations, which are used in this paper,  we refer the reader to
\cite{SW} and \cite{SKM}.

\section{Background Material}\label{secback}

\subsection{Iteration Infinitely Divisible and Iteration Stable Random Tessellations}\label{secstit}

Consider a compact and convex set $W\subset{\Bbb R}^2$ with non-empty interior and a diffuse
(non-atomic) measure $\Lambda$ on the space of $[{\Bbb R}^2]$ of lines in the plane enjoying the
{\it local finiteness property}, stating that $\Lambda([C])<\infty$ for any compact 
$C\subset{\Bbb R}^2$, under the usual notation
$$[C]=\{\text{lines}\ \ L:\ \ L\cap C\neq\emptyset\},$$
letting $[C]$ stand for the set of lines that have non-empty intersection with $C$.
Below we briefly describe the construction of an iteration
infinitely divisible tessellation in $W,$ the construction is called the MNW-construction to
honor its inventors \--- Mecke, Nagel and Wei\ss\ , who introduced it in \cite{NW05}. 
To begin, assign to $W$ an exponentially distributed random lifetime with parameter given by 
$\Lambda([W]).$ Upon expiry of this random time
the cell $W$ dies and splits into two sub-cells $W^+$ and $W^-$ separated by a line in $[W]$
chosen according to the law $\Lambda([W])^{-1}\Lambda(\cdot\cap[W])$. The resulting new
cells $W^+$ and $W^-$ are again assigned independent exponential lifetimes with respective
parameters $\Lambda([W^+])$ and $\Lambda([W^-])$, whereupon the construction continues
recursively and independently in each of the sub-cells $W^+$ and $W^-$ until some deterministic
time threshold $t>0$ is reached. The random tessellation constructed in $W$ until the time $t$ will
be denoted by $Y(t,W)$. To put it formally, by $Y(t,W)$ we understand here the random closed subset
of $W$ arising as the union of the boundaries of cells constructed by the  time $t.$ The cell-splitting edges are called the \textit{maximal edges} (in the related literature often called \textit{I-segments}, as assuming shapes similar to that of the literal I) of the tessellation $Y(t,W)$ and the family of all such edges is denoted by $\MaxEdges(Y(t,W))$. Note that such maximal edges can get further subdivided between their birth time and the time $t$, that is to say in the course of the MNW dynamics there can appear additional vertices in their relative interiors. Thus, we have a distinction between maximal edges with a possible interior structure and those edges which are not maximal (for example the sides of cells or the primitive elements that are bounded by vertices but have no interior structure).\\ In general, the random tessellations $Y(t,W)$ in $W$ do not have to arise as windowed restrictions of stationary (stochastically translation invariant) or isotropic (stochastically rotation invariant) processes. If we assume in addition though that the driving measure $\Lambda$ has, under polar parametrization, the product structure $\Lambda=\tau \ell_+\otimes{\cal R}$ with $\tau$ a positive constant, $\ell_+$ the Lebesgue measure on ${\Bbb R}_+$ and a spherical directional distribution $\cal R$ on the unit circle ${\cal S}_1$, then $Y(t,W)$ has the property of being a suitable windowed restriction of a stationary random field. The random tessellation $Y(t,W)$ is furthermore isotropic in a similar restriction-wise sense if and only if $\cal R$ is the uniform distribution, i.e. if and only if $\Lambda$ is a multiple of the standard motion-invariant measure $\Lambda_{iso} := \ell_+ \otimes
\nu_1$ of unit length intensity on the space of lines in the plane, where $\nu_1$ is the uniform distribution on ${\cal S}_1.$\\ We review now some of the important properties of the tessellations $Y(t,W)$, proofs of which can be found in \cite{MNW}, \cite{NW05} and \cite{ST}. We formulate them only for the planar case even if they are true in higher dimensions as well.
\begin{enumerate}
  \item The random tessellations $Y(t,W)$ are consistent in $W,$ that is to say
                  $Y(t,V) \cap W \overset{D}{=} Y(t,W)$ for $W \subseteq V$ for $V,W$ compact convex, where $\overset{D}{=}$ means equality in distribution. This implies
                  the existence of the whole-plane tessellation $Y(t)$ such that $Y(t,W) \overset{D}{=} Y(t) \cap W$
                  for each compact convex $W.$ 
  \item The MNW-construction satisfies the Markov property with respect to \textit{iteration} of tessellations in the time parameter $t$, i.e $$Y(s,W)\boxplus Y(t,W)\overset{D}{=}Y(s+t,W),$$ where $\boxplus$ denotes the operation of iteration for whose exact definition we refer to the above mentioned papers. This operation can roughly be explained as follows: Let $Y_0$ be a random tessellation and let $T_1,T_2,\ldots$ be a sequence of independent and identically distributed random tessellations in the plane. To each cell $c_k$ of the \textit{frame tessellation} $Y_0$ we associate the tessellation $T_k$. Now, we make a local superposition of $Y_0$ and the cutouts of $T_k$ in each cell $c_k$ of the frame tessellation. If we let $Y_1$ be a tessellation with the same distribution as $T_k$ for any $k=1,2,\ldots$, we denote the outcome of this procedure of local superposition by $Y_0\boxplus Y_1$ and call it the \textit{iteration} of $Y_0$ with $Y_1$. 
 \item It directly follows from the above Markov property that the random tessellations $Y(t,W)$ are infinitely divisible with respect to iteration, i.e. $$Y(t,W)\overset{D}{=}\underbrace{Y(t/n,W)\boxplus\ldots\boxplus Y(t/n,W)}_{n}$$ for any $n\in{\Bbb N}$ and any compact convex $W$. For this reason, the tessellations $Y(t,W)$ are referred to as \textit{iteration infinitely divisible random tessellations}. However, it is worth pointing out that it is currently not known whether \textit{any} iteration infinitely divisible random tessellation can be obtained by the MNW-construction. If in addition the driving measure $\Lambda$ is translation invariant, it can be verified that $Y(t,W)$ is even stable under iteration, i.e. $$Y(t)\overset{D}{=}n(\underbrace{Y(t)\boxplus\ldots\boxplus Y(t)}_{n}),\ \ \ \ \ n\in{\Bbb N},$$ where $n(\cdot)$ means the dilation with a factor $n$, i.e. $nY=\{n\cdot x: x\in Y\}$. Note that this equation must be understood symbolically, since properly the operation of iteration involves a sequence of independent and identically distributed random tessellations, but we adopt here the usual convention from the earlier work \cite{MNW}. Thus, in the stationary case we are in fact dealing with \textit{iteration stable random tessellations} or \textit{STIT-tessellations} for short.
	\item The intersection of an iteration infinitely divisible random tessellation $Y(t,W)$ having driving measure $\Lambda$ with an arbitrary line $L\in[W]$ is a Poisson point process with intensity measure $A\mapsto t\Lambda([A])$ with $A\subset W\cap L$ a Borel set. In particular,  for $x,y\in W$ the probability that they belong to the same cell of $Y(t,W)$ equals $${\Bbb P}(\text{$x$ and $y$ are in the same cell of $Y(t,W)$})=e^{-t\Lambda([\overline{xy}])},$$ where by $\overline{xy}\subset W$ we mean the line segment joining the points $x$ and $y$.
\end{enumerate}

\subsection{Martingales in the MNW Construction}\label{secmart}

As already seen in the introduction above, the MNW-construction of iteration infinitely divisible or iteration stable random tessellations $Y(t,W)$ in finite volumes $W\subset{\Bbb R}^2$ for general locally finite and diffuse
driving measures $\Lambda$ enjoys the Markov property in the continuous time parameter $t.$ 
In our previous work \cite{ST} we
have used this fact combined with the classical theory of martingale problems for pure jump Markov processes
to construct a class of natural martingales associated to the MNW process. In this paper we only need
a part of that theory. To formulate the required results, for a bounded and measurable functional $\phi$
of a line segment (tessellation edge) and for a tessellation $Y,$ usually taken to be a realization
of $Y(t,W)$ for some $t > 0,$ write 
\begin{equation}\nonumber \Sigma_{\phi}(Y)=
 \sum_{e\in\MaxEdges(Y)}\phi(e).
\label{defsigma}
\end{equation}
Note, that with $Y$ as above, for each straight line $L \in [W]$ the intersection $Y \cap L$ is just
a one-dimensional {\it tessellation} of $L \cap W$ which can be identified with the collection
$\Segments(Y \cap L)$ of its constituent segments. Bearing this in mind we write
\begin{equation}
 A_{\phi}(Y) = \int_{[W]} \sum_{e \in \Segments(Y \cap L \cap W)} \phi(e) \Lambda(dL).
\label{defA}
\end{equation}
It is also convenient to introduce the bar notation for centered versions of these quantities with
$Y = Y(t,W),$ that is to say $$\bar{\Sigma}_{\phi}(Y(t,W)) := \Sigma_{\phi}(Y(t,W)) -
 {\Bbb E}\Sigma_{\phi}(Y(t,W))$$
and likewise $$\bar{A}_{\phi}(Y(t,W)) := A_{\phi}(Y(t,W)) - {\Bbb E}A_{\phi}(Y(t,W)).$$
With this notation, in view of (41) in \cite{ST}  we have
\begin{proposition}\label{MartProp}
 For bounded measurable segment functionals $\phi$ and $\psi$
 the stochastic processes
 \begin{equation}\label{ORD1mart}
  \Sigma_{\phi}(Y(t,W)) - \int_0^t A_{\phi}(Y(s,W)) ds 
 \end{equation}
  and
\begin{eqnarray}
\nonumber  & & \bar\Sigma_{\phi}(Y(t,W)) \bar\Sigma_{\psi}(Y(t,W)) - \int_0^t A_{\phi \psi}(Y(s,W)) ds -\\
 &-& \int_0^t [\bar{A}_{\phi}(Y(s,W)) \bar\Sigma_{\psi}(Y(s,W)) + \bar{A}_{\psi}(Y(s,W))\bar\Sigma_{\phi}(Y(s,W))]ds\label{VARMART}
\end{eqnarray}
are martingales with respect to the filtration $\Im_t$ generated by $(Y(s,W))_{0 \leq s \leq t}.$
\end{proposition}

\subsection{Mean Values for Edge and Vertex Processes}\label{secmv}

This subsection recapitulates some basic first order properties of iteration infinitely divisible tessellations,
mostly known in the stationary set-up, for sake of reference in further sections.  Using (\ref{ORD1mart}) for  $\phi \equiv 1$ yields, upon taking expectations,
 \begin{equation}\label{ES1}
  {\Bbb E}\Sigma_1(Y(t,W)) = \int_0^t {\Bbb E}A_1(Y(s,W)) ds. 
 \end{equation}
 However, in view of the definition (\ref{defA}) it is easily verified that 
 \begin{equation}\label{A1}
   A_1(Y) = \Lambda([W])+\sum_{e \in \MaxEdges(Y)} \Lambda([e]) = \Lambda([W])+\Sigma_{\Lambda[\cdot]}(Y),
 \end{equation}
 where $\Lambda[\cdot]$ stands for the edge functional $e \mapsto \Lambda([e]).$ 
 Applying (\ref{ORD1mart}) once again for $\phi = \Lambda[\cdot]$ leads to
 $$ 
  {\Bbb E}\Sigma_{\Lambda[\cdot]}(Y(t,W)) = \int_0^t {\Bbb E}A_{\Lambda[\cdot]}(Y(s,W)) ds. $$
But 
\begin{equation}\label{ACONST}
 A_{\Lambda[\cdot]}(Y(s,W)) = \int_{[W]} \Lambda([L \cap W]) \Lambda(dL) = \lparen
 \Lambda \cap \Lambda \rparen(W),\ \ \ \ \forall s\in(0,t],
\end{equation}
where the locally finite {\it point-intersection measure} $\lparen \Lambda \cap \Lambda \rparen$
on ${\Bbb R}^2$ is given by
\begin{equation}\label{IntersectMeas}
 \lparen \Lambda \cap \Lambda \rparen := \int_{[{\Bbb R}^2]} \int_{[L]} \delta_{L \cap L'} \Lambda(dL'),
 \Lambda(dL)
 \end{equation}
 where, recall, $\delta_{(\cdot)}$ stands for the Dirac unit mass at the argument, so that in other words
 $$ \lparen \Lambda \cap \Lambda \rparen(A) = (\Lambda \times \Lambda)\{ (L_1,L_2) \in [A]\times[A],\;
  L_1 \cap L_2 \in A \},\ \ \ \ A \subseteq {\Bbb R}^2. $$
 Hence
 \begin{equation}\label{ESLA}
  {\Bbb E}\Sigma_{\Lambda[\cdot]}(Y(s,W)) = s \lparen \Lambda \cap \Lambda \rparen(W).
 \end{equation}
 Combining (\ref{ES1}), (\ref{A1}) and (\ref{ESLA}) finally yields
 \begin{equation}\label{ES1fin}
  {\Bbb E}\Sigma_{1}(Y(t,W)) = t\Lambda([W])+\frac{t^2}{2} \lparen \Lambda \cap \Lambda \rparen(W).
 \end{equation}
Note, that the formula (\ref{ES1fin}) when specialized to the translation-invariant set-up contains an
extra boundary correction term $t \Lambda([W])$ in comparison to the classical mean value formula
for the iteration stable (STIT) random tessellations as given in \cite{NW06}, which says that the density of
maximal edges in $W$ is just $\frac{t^2}{2} \lparen \Lambda \cap \Lambda \rparen(W)$. This additional
boundary correction term $t\Lambda([W])$ comes from the fact that we count edges rather than edge
midpoints. Thus, it can happen that in two neighboring regions one observes two distinct edges which
may coalesce into one edge when putting these regions together into one area.

\section{Second Order Theory for Edge and Vertex Processes}\label{secvar}

\subsection{Variance Calculation for the General Case}\label{secvargen}

 We consider the most general case first and study iteration infinitely divisible random tessellations
 $Y(t,W)$ with general locally finite and non-atomic driving measures $\Lambda.$ We fix $t>0$
 and a compact and convex  observation window $W\subset{\Bbb R}^2$ as in Section \ref{secstit}.\\ 
 First we use (\ref{VARMART}) with $\phi=\psi\equiv 1$ to conclude that
 $$ \bar\Sigma_1^2(Y(t,W))-\int_0^t A_1(Y(s,W))ds-2\int_0^t\bar{A}_1(Y(s,W))\bar\Sigma_1(Y(s,W))ds $$
 is a martingale with respect to $\Im_t$. Using (\ref{A1}) and taking expectations we get
\begin{eqnarray}
\nonumber \Var(\Sigma_1(Y(t,W))) &=& t\Lambda([W])+\int_0^t {\Bbb E} \Sigma_{\Lambda[\cdot ]}(Y(s,W)) ds\\
& & +2 \int_0^t \Cov(\Sigma_{\Lambda[\cdot]}(Y(s,W)),\Sigma_1(Y(s,W))) ds. \label{VARFOR1}
\end{eqnarray}
It remains to find an expression for the covariance $\Cov(\Sigma_{\Lambda[\cdot]}(Y(s,W)),\Sigma_1(Y(s,W))).$ Such an expression can be found by applying (\ref{VARMART}) once again, this time with $\phi=1$ and $\psi=\Lambda[\cdot]$ and $t$ replaced by $s$ and $s$ by $u$. We obtain in this way
\begin{eqnarray}
\nonumber \Cov(\Sigma_{\Lambda[\cdot]}(Y(s,W)),\Sigma_1(Y(s,W))) &=& \int_0^s {\Bbb E} A_{\Lambda[\cdot]}(Y(u,W)) du\\
\nonumber &+& \int_0^s \Cov(A_1(Y(u,W)),\Sigma_{\Lambda[\cdot]}(Y(u,W))) du\\
\nonumber &+& \int_0^s \Cov(A_{\Lambda[\cdot]}(Y(u,W)),\Sigma_1(Y(u,W))) du.
\end{eqnarray}
In view of (\ref{ACONST}), $A_{\Lambda[\cdot]}(\cdot)$ is a constant
and hence the covariance involving it vanishes. Resorting again to (\ref{A1}) we end up with
\begin{eqnarray}
\nonumber & & \Cov(\Sigma_{\Lambda[\cdot]}(Y(s,W)),\Sigma_1(Y(s,W)))\\
&=& \int_0^s {\Bbb E} A_{\Lambda[\cdot]}(Y(u,W)) du + \int_0^s \Var(\Sigma_{\Lambda[\cdot]}(Y(u,W))) du.\label{COVFOR2}
\end{eqnarray}
Putting together (\ref{VARFOR1}) with (\ref{COVFOR2}) yields the following expression
for $\Var(\Sigma_1(Y(t,W)))$:
\begin{eqnarray}
\nonumber & & \Var(\Sigma_1(Y(t,W)))= t \Lambda([W]) + \int_0^t{\Bbb E}\Sigma_{\Lambda[\cdot]}(Y(s,W))ds\\
&+& 2\left(\int_0^t\int_0^s {\Bbb E}A_{\Lambda[\cdot]}(Y(u,W))duds+\int_0^t\int_0^s\Var(\Sigma_{\Lambda[\cdot]}(Y(u,W))duds\right).\label{VARSI1}
\end{eqnarray}
 It remains to find $\Var(\Sigma_{\Lambda[\cdot]}(Y(t,W))).$ To find an expression, apply (\ref{VARMART})
 again with $\phi = \psi = \Lambda[\cdot]$ to get, upon taking expectations and using 
 that $A_{\Lambda[\cdot]}$ is a constant as remarked in (\ref{ACONST}) above, 
 \begin{equation}\label{VARSL}
  \Var(\Sigma_{\Lambda[\cdot]}(Y(t,W))) = \int_0^t {\Bbb E} A_{(\Lambda[\cdot])^2}(Y(s,W)) ds. 
 \end{equation}
  However, using (\ref{defA}) 
  $$ A_{(\Lambda[\cdot])^2}(Y(s,W)) = \int_{[W]} \int_{[L \cap W]} \int_{[L \cap W]} 
     {\bf 1}[L \cap L_1 \mbox{ and }
          L \cap L_2 \mbox{ are } $$ $$ \mbox{\hspace{5cm} in the same cell of } Y(s,W)] \Lambda(dL_1) \Lambda(dL_2) \Lambda(dL) $$
 and hence
\begin{eqnarray} 
\nonumber & & {\Bbb E}A_{(\Lambda[\cdot])^2}(Y(s,W))\\
\nonumber &=& \int_{[W]} \int_{[L \cap W]} \int_{[L \cap W]} {\Bbb P}(L \cap L_1 \mbox{ and }
          L \cap L_2 \mbox{ are}\\
\nonumber & & \mbox{\hspace{3.2cm} in the same cell of } Y(s,W)) \Lambda(dL_1) \Lambda(dL_2) \Lambda(dL)\\
  &=& \int_{[W]} \int_{[L \cap W]} \int_{[L \cap W]} \exp\left(-s\Lambda\left(\left[ L(L_1,L_2)
  \right]\right)\right) 
  \Lambda(dL_1) \Lambda(dL_2) \Lambda(dL)\label{ALASQ}
\end{eqnarray}
 where $L(L_1,L_2)$ stands for the segment joining the points $L\cap L_1$ and $L\cap L_2$ and where the last equality follows by the property 4 of the tessellation $Y(t,W)$ as listed in Section \ref{secstit}.\\ To neatly formulate our theory, denote by
  $\lparen (\Lambda \times \Lambda) \cap \Lambda \rparen$ the {\it segment-intersection 
  measure} on the space $\bar{[} {\Bbb R}^2 \bar{]}$ of finite linear segments in ${\Bbb R}^2$
  given by
  \begin{equation}\label{SIntersectMeas}
   \lparen (\Lambda \times \Lambda) \cap \Lambda \rparen = \int_{[{\Bbb R}^2]} \int_{[L]} \int_{[L]}
   \delta_{L(L_1,L_2)} \Lambda(dL_1) \Lambda(dL_2) \Lambda(dL)
  \end{equation}
  and observe that this defines a locally finite measure, charging finite mass on collections $\bar{[} A \bar{]}$
  of segments  with both ends falling into a bounded set $A \subset {\Bbb R}^2,$ because of the local
  finiteness of $\Lambda.$
  With this notation, combining (\ref{ACONST}), (\ref{ESLA}),  (\ref{COVFOR2}), (\ref{VARSI1}), (\ref{VARSL}) and
  (\ref{ALASQ}) yields the following result:
 \begin{theorem}\label{thmvar}
  For general locally finite and diffuse driving measures $\Lambda,$ denoting by
  $$ T_n^{\exp}(u) = \sum_{k=n}^{\infty} u^k/k! = \exp(u) - \sum_{k=0}^{n-1} u^k/k! $$
  the $n$-th tail of the exponential series at $u,$ we have
  \begin{equation}\nonumber\label{VARSLthm}
   \Var(\Sigma_{\Lambda[\cdot]}(Y(t,W))) = -\int_{\bar{[} W \bar{]}}
    \frac{T_1^{\exp}(-t\Lambda([e]))}{\Lambda([e])} 
   \lparen (\Lambda \times \Lambda) \cap \Lambda \rparen(de)
  \end{equation}
  and
  \begin{eqnarray}
   \nonumber & & \Cov(\Sigma_{\Lambda[\cdot]}(Y(t,W)),\Sigma_1(Y(t,W)))\\
   \nonumber &=& t \lparen \Lambda \cap \Lambda \rparen(W)+  
   \int_{\bar{[} W \bar{]}} 
   \frac{T^{\exp}_2(-t\Lambda([e]))}{\Lambda([e])^2}
   \lparen (\Lambda \times \Lambda) \cap \Lambda \rparen(de)\label{COVSLS1thm}
  \end{eqnarray}
  and
  \begin{eqnarray}
  \nonumber \Var(\Sigma_1(Y(t,W))) &=& t \Lambda([W]) +  \frac{3 t^2}{2}
   \lparen \Lambda \cap \Lambda \rparen(W)\\
   & & - 2 \int_{\bar{[}W\bar{]}}
   \frac{T^{\exp}_3(-t\Lambda([e]))}{\Lambda([e])^3}
   \lparen (\Lambda \times \Lambda) \cap \Lambda \rparen(de).\label{VARS1thm} 
 \end{eqnarray}
 \end{theorem}
 
\subsection{Vertex Pair-Correlations for the General Case}\label{secgenPCF}

 Also in this subsection we stay in the general set-up of locally finite and diffuse $\Lambda.$ 
 We shall extend here the calculations made in Subsection \ref{secvargen} above to determine
 the pair-correlation structure of the vertex point process ${\cal V}_{Y(t,W)}$ generated
 by $Y(t,W).$ For definiteness we adopt the convention that ${\cal V}_{Y(t,W)}$ does
 not include the boundary vertices, this way each vertex arises at intersection of exactly two
 maximal edges. Recalling the consistency relation $Y(t,W) = Y(t) \cap W$, we see that the
 covariance structure between bounded regions $U,V \subset {\Bbb R}^2$ does not
 depend on $W$ as soon as both $U$ and $V$ are contained in the interior of $W.$
 To put this in formal terms, consider the whole-plane covariance measure $\Cov({\cal V}_{Y(t)})$
 of the point process ${\cal V}_{Y(t)}$ on $({\Bbb R}^2)^2 = {\Bbb R}^2 \times {\Bbb R}^2$ (also called the second-order cumulant measure)
 given by the relation
 \begin{equation}\label{CovMeas}
  \int_{({\Bbb R}^2)^2} (f \otimes g) d \Cov({\cal V}_{Y(t)}) =
  \Cov(\Sigma_{\eta^{f}}(Y(t)),\Sigma_{\eta^{g}}(Y(t))) 
 \end{equation}
 holding for all $f,g : {\Bbb R}^2 \to {\Bbb R}$ bounded measurable and of bounded support,
 where $\eta^f$ is the edge functional
 $$ \eta^f(e) = \sum_{x \in \Vertices(e)} f(x) $$
 and likewise for $\eta^g.$ Note, that, even though we are apparently dealing with functionals
 $\Sigma_{\eta^{f}}(Y(t))$ and $\Sigma_{\eta^{g}}(Y(t))$ defined on the whole-plane process,
 they can be safely replaced by  $\Sigma_{\eta^{f}}(Y(t,W))$ and $\Sigma_{\eta^{g}}(Y(t,W))$,
 respectively, for some $W$ containing the supports of $f$ and $g,$ hence our martingale 
 relations given in Proposition \ref{MartProp} hold here with no extra assumptions. For the
 same reasons all integrals below with apparently unbounded integration domains are
 effectively bounded due to the bounded supports of $f$ and $g,$ which we are going to
 exploit without further mention. It is readily seen from (\ref{defA}) that for each possible
 realization $Y$ of $Y(t)$ or $Y(t,W)$ in a domain $W$ containing the supports of $f,$
 \begin{equation}\label{AETAF}
  A_{\eta^f}(Y) = 2 \int_{[{\Bbb R}^2]} \sum_{x \in L \cap Y} f(x) \Lambda(dL) =
  2 \Sigma_{\Lambda^f[\cdot]}(Y),
 \end{equation}
 where $\Lambda^f[e] = \int_{[e]} f(e \cap L) \Lambda(dL)$ and the factor $2$ comes from the fact that each point of the tessellation is contained in exactly two maximal edges. 
 Consequently, using (\ref{VARMART}) for $\phi = \eta^f, \;\psi = \eta^g$ and taking expectations, we get
 \begin{eqnarray}
 \nonumber \Cov(\Sigma_{\eta^{f}}(Y(t)),\Sigma_{\eta^{g}}(Y(t))) &=& \int_0^t {\Bbb E}A_{\eta^f \eta^g}(Y(s)) ds \\
 \nonumber & & + 2\int_0^t \Cov(\Sigma_{\Lambda^f[\cdot]}(Y(s)),\Sigma_{\eta^g}(Y(s))) ds\\
  & & + 2\int_0^t \Cov(\Sigma_{\Lambda^g[\cdot]}(Y(s)),\Sigma_{\eta^f}(Y(s))) ds.\label{PCcovar}
 \end{eqnarray}
 Proceeding as in the previous Subsection \ref{secvargen}, we turn now to the calculation of  
 the covariance $\Cov(\Sigma_{\Lambda^f[\cdot]}(Y(s)),\Sigma_{\eta^g}(Y(s))).$
 To this end we note that
 \begin{equation}\label{ALAF}
  A_{\Lambda^f[\cdot]}(Y) = \int_{[{\Bbb R}^2]} \Lambda^f[L] \Lambda(dL) = \int_{{\Bbb R}^2}
  f d\lparen \Lambda \cap \Lambda \rparen, 
 \end{equation}
 whence $A_{\Lambda^f[\cdot]}(\cdot)$ is a constant, 
 and we use again (\ref{VARMART}) for $\phi = \Lambda^f,\; \psi = \eta^g$ to get in view of (\ref{AETAF})
 \begin{eqnarray}
 \nonumber & & \Cov(\Sigma_{\Lambda^f[\cdot]}(Y(s)),\Sigma_{\eta^g}(Y(s)))\\
 &=& \int_0^s {\Bbb E} A_{\Lambda^f[\cdot] \eta^g}(Y(u)) du 
 + 2\int_0^s \Cov(\Sigma_{\Lambda^f[\cdot]}(Y(u)),\Sigma_{\Lambda^g[\cdot]}(Y(u))) du.\label{PCcovar2}
 \end{eqnarray}
 Finally, one further use of (\ref{VARMART}) with $\phi = \Lambda^f,\; \psi = \Lambda^g$ and application
 of (\ref{ALAF}) yields
 \begin{equation}\label{PCcovar3}
  \Cov(\Sigma_{\Lambda^f[\cdot]}(Y(u)),\Sigma_{\Lambda^g[\cdot]}(Y(u))) =
  \int_0^u {\Bbb E}A_{\Lambda^f[\cdot] \Lambda^g[\cdot]}(Y(v)) dv.
 \end{equation}
 Thus, using (\ref{PCcovar}), then twice (\ref{PCcovar2}), once with $f$ and $g$ interchanged, and then
 (\ref{PCcovar3}), we get
 $$ \Cov(\Sigma_{\eta^{f}}(Y(t)),\Sigma_{\eta^{g}}(Y(t))) = \int_0^t {\Bbb E} A_{\eta^f \eta^g}(Y(s)) ds$$
 \begin{equation}
 +2 \int_0^t \int_0^s {\Bbb E}A_{\Lambda^f[\cdot] \eta^g + \Lambda^g[\cdot] \eta^f}(Y(u)) du
 +8  \int_0^t \int_0^s \int_0^u {\Bbb E} A_{\Lambda^f[\cdot]\Lambda^g[\cdot]}(Y(v)) dv du ds,\label{PCcovar4}
 \end{equation}
 since $A_{\Lambda^f[\cdot]\eta^g}+A_{\Lambda^g[\cdot]\eta^f}=A_{\Lambda^f[\cdot]\eta^g+\Lambda^g[\cdot]\eta^f}$.
 It remains to calculate the expectations of the $A_{(\cdot)}$ functionals present in these integrals.
 However, this is easily done by recalling that, for $L \in [{\Bbb R}^2],$ the intersection $Y(t) \cap L$
 is the Poisson point process with intensity measure $L \supseteq A \mapsto  t\Lambda([A]),$
 see property 4 in Subsection \ref{secstit}, whence, in view of (\ref{defA}),
 \begin{eqnarray}
  \nonumber {\Bbb E}A_{\phi}(Y(t)) &=& \frac{1}{2}\int_{\bar{[} {\Bbb R}^2 \bar{]}} \phi(e) \exp(-t\Lambda([e])) \lparen (t\Lambda \times t\Lambda) \cap \Lambda\rparen(de)\\
  &=&\frac{t^2}{2} 
  \int_{\bar{[} {\Bbb R}^2 \bar{]}} \phi(e) \exp(-t\Lambda([e])) \lparen (\Lambda \times
  \Lambda) \cap \Lambda \rparen(de)\label{GenAExp}
 \end{eqnarray}
 for bounded measurable $\phi$ such that $\phi(\emptyset) = 0$ and locally defined in the sense that there
 exists a bounded convex $W$ such that $\phi(e) = \phi(e \cap W)$ for all $e.$ Note that 
 the extra prefactor of $1/2$ comes from the fact that the segment-intersection measure defined by (\ref{SIntersectMeas}) counts each segment twice, once for each of the 
 two orderings of its two termini.
 Putting (\ref{GenAExp}) together with (\ref{PCcovar4}) yields now
 \begin{eqnarray}
 \nonumber & & \Cov(\Sigma_{\eta^{f}}(Y(t)),\Sigma_{\eta^{g}}(Y(t))) =\\
 \nonumber &+& {1\over 2} \int_{\bar{[} {\Bbb R}^2 \bar{]}}
      \eta^f(e) \eta^g(e) \ {\cal I}^1(s^2 \exp(-t\Lambda([e]));t) \ \lparen (\Lambda \times
      \Lambda) \cap \Lambda \rparen(de)\\
 \nonumber &+& \int_{\bar{[} {\Bbb R}^2 \bar{]}} (\Lambda^f[e] \eta^g[e] + \Lambda^g[e] \eta^f(e)) \
        {\cal I}^2(s^2 \exp(-s\Lambda([e]));t) \ \lparen (\Lambda \times
      \Lambda) \cap \Lambda \rparen(de)\\
 &+& 4 \int_{\bar{[} {\Bbb R}^2 \bar{]}} \Lambda^f[e] \Lambda^g[e] \ {\cal I}^3(s^2 \exp(-s\Lambda([e]));t)
  \ \lparen (\Lambda \times \Lambda) \cap \Lambda \rparen(de), \label{PCcovar5}
 \end{eqnarray}
 where the multiple integral ${\cal I}^n$ is given by
 $$ {\cal I}^n(f(s);t) := \int_0^t  \int_0^{s_1} \ldots \int_0^{s_{n-1}} f(s) ds ds_{n-1} \ldots ds_1 =
     \frac{1}{(n-1)!} \int_0^t (t-s)^{n-1} f(s) ds $$
 so that in particular 
 \begin{eqnarray}
 \label{I1} {\cal I}^1(s^2 \exp(-\lambda s);t) &=& \lambda^{-3} (2 - (\lambda^2 t^2 + 2 \lambda t + 2) \exp(-\lambda t)),\\
 \label{I2} {\cal I}^2(s^2 \exp(-\lambda s);t) &=& \lambda^{-4}(2\lambda t - 6 + (\lambda^2 t^2 + 4\lambda t + 6) \exp(-\lambda t)),\\
 \label{I3} {\cal I}^3(s^2 \exp(-\lambda s);t) &=& \lambda^{-5}(\lambda^2 t^2 - 6 \lambda t +12 - (\lambda^2 t^2 + 6 \lambda t + 12) \exp(-\lambda t)).
 \end{eqnarray}
 For a segment (edge) $e$ consider the measures $\Delta^e$ and $\Lambda[\cdot \cap e]$ on ${\Bbb R}^2$ that are defined by
 $$ \Delta^e := \sum_{x \in \Vertices(e)} \delta_x $$
 and
 $$ (\Lambda[\cdot \cap e])(A) = \Lambda([A \cap e]),\ \ \ A \subseteq {\Bbb R}^2. $$  
 With this notation, putting together (\ref{CovMeas}) and (\ref{PCcovar5})
 yields in view of the definitions of $\eta^f$ and $\Lambda^f[\cdot]$ the following
 \begin{theorem}\label{CovMeasThm}
  For general locally finite and diffuse driving measures $\Lambda$ we have
  \begin{eqnarray}
 \nonumber & &
  \Cov({\cal V}_{Y(t)}) =  \int_{\bar{[} {\Bbb R}^2 \bar{]}}
   \frac{1}{2} (\Delta^e \otimes \Delta^e)  {\cal I}^1(s^2 \exp(-s\Lambda([e]));t)
   \lparen (\Lambda \times \Lambda) \cap \Lambda \rparen(de) + \\
   \nonumber & &
      \int_{\bar{[} {\Bbb R}^2 \bar{]}}  
         (\Delta^e \otimes \Lambda[\cdot \cap e] + \Lambda[\cdot \cap e] \otimes \Delta^e) \
         {\cal I}^2(s^2 \exp(-s\Lambda([e]));t) \lparen (\Lambda \times \Lambda) \cap \Lambda \rparen(de) + \\
      & & \ \ \ \ \ \  4 \int_{\bar{[} {\Bbb R}^2 \bar{]}}  
             (\Lambda[\cdot \cap e] \otimes \Lambda[\cdot \cap e])
             {\cal I}^3(s^2 \exp(-s\Lambda([e]));t) \lparen (\Lambda \times \Lambda) \cap \Lambda \rparen(de).
\label{CovMeasExpr}
  \end{eqnarray}
\end{theorem}
 An intuitive understanding of the structure of the covariance measure in Theorem \ref{CovMeasThm}
 comes by noting that the first integral in (\ref{CovMeasExpr}) takes into account 
 pairs of vertices constituting ends of the same maximal edge, the second one corresponds to pairs 
 of vertices with the property that one of them is an internal vertex of a maximal edge of which the second
 point is a terminus, whereas the third integral corresponds to pairs of vertices lying on the
 same maximal edge but not being its termini. Thus, other pairs of vertices (not lying on the
 same maximal edge) are not present in the covariance measure, roughly speaking this is because
 the maximal edges are {\it the only means of propagating dependencies} in iteration infinitely divisible
 tessellations, an intuition to be made more specific in our further work in progress.

\subsection{Edge-Vertex Correlations in the General Case}\label{secgenEVPCF}

 In this subsection, still placing ourselves in the general setting of a locally finite and diffuse $\Lambda,$
 we consider the covariance measure between the vertex point process and edge length process
 generated by $Y(t).$ To this end, define the (random) edge-length measure ${\cal E}_{Y(t)}$
 of $Y(t)$ by putting for bounded Borel $A \subseteq {\Bbb R}^2$
 $$ {\cal E}_{Y(t)}(A) = \sum_{e \in \MaxEdges(Y(t))} \ell(e \cap A) $$
 with $\ell(\cdot)$ standing for the usual one-dimensional length. The object of our
 interest is the measure $\Cov({\cal V}_{Y(t)},{\cal E}_{Y(t)})$ given by
 \begin{equation}\label{CovMixed}
  \int_{({\Bbb R}^2)^2} (f \otimes g) d \Cov({\cal V}_{Y(t)},{\cal E}_{Y(t)}) = 
  \Cov(\Sigma_{\eta^f}(Y(t)),\Sigma_{J^g}(Y(t)))
 \end{equation}
 for bounded measurable $f,g : {\Bbb R}^2 \to {\Bbb R}$ with bounded support, where
 $J^g$ denotes the functional $ J^g(e) = \int_e g(x) \ell(dx).$ 
 Similarly as in (\ref{ALAF}) we have
 \begin{equation}\label{AIF}
   A_{J^g}(Y) = \int_{[{\Bbb R}^2]}  J^g(L) \Lambda(dL).
 \end{equation}
 Thus, $A_{J^g}$ is constant and hence, using (\ref{VARMART}) for $\phi = \eta^f$ and $\psi = J^g,$
 taking expectations and recalling (\ref{AETAF}), yields
 $$\Cov(\Sigma_{\eta^f}(Y(t)),\Sigma_{J^g}(Y(t))) = \int_0^t 
 {\Bbb E} A_{\eta^f J^g}(Y(s)) ds + $$
\begin{equation}\label{CovMix1}
  2  \int_0^t \Cov(\Sigma_{\Lambda^f[\cdot]}(Y(s)),\Sigma_{J^g}(Y(s))) ds.
 \end{equation}
 Using (\ref{VARMART}) once again, with $\phi = \Lambda^f[\cdot]$ and $\psi = J^g,$ upon
 taking expectations and recalling (\ref{AIF}) and (\ref{ALAF}), we get
 \begin{equation}\nonumber\label{CovMix2}
   \Cov(\Sigma_{\Lambda^f[\cdot]}(Y(s)),\Sigma_{J^g}(Y(s))) =
   \int_0^s {\Bbb E} A_{\Lambda^f[\cdot] J^g}(Y(u)) du. 
 \end{equation}
 Substituting into (\ref{CovMix1}) leads us to
 \begin{eqnarray}
  \nonumber \Cov(\Sigma_{\eta^f}(Y(t)),\Sigma_{J^g}(Y(t))) &=&  \int_0^t 
  {\Bbb E} A_{\eta^f J^g}(Y(s)) ds\\
  \nonumber & & + 2 \int_0^t \int_0^s {\Bbb E} A_{\Lambda^f[\cdot] J^g}(Y(u)) duds. \label{CovMix3}
 \end{eqnarray}
 Applying (\ref{GenAExp}), we obtain therefore 
\begin{eqnarray} 
\nonumber & & \Cov(\Sigma_{\eta^f}(Y(t)),\Sigma_{J^g}(Y(t)))\\  
\nonumber & = & \frac{1}{2} \int_{\bar{[} {\Bbb R}^2 \bar{]}} \eta^f(e) J^g(e) \ {\cal I}^1(s^2 \exp(-s\Lambda([e]));t)
       \lparen (\Lambda \times \Lambda) \cap \Lambda \rparen(de)\\
  && + \int_{\bar{[} {\Bbb R}^2 \bar{]}} \Lambda^f[e] J^g(e) \ {\cal I}^2(s^2 \exp(-s\Lambda([e]));t)
       \lparen (\Lambda \times \Lambda) \cap \Lambda \rparen(de).\label{CovMix4}
 \end{eqnarray}
 Consequently, putting (\ref{CovMix4}) together with (\ref{CovMixed}) and defining the measure 
 $ (\ell(\cdot \cap e))(A) := \ell(A \cap e),\; A \subseteq {\Bbb R}^2, $
 we obtain
 \begin{theorem}\label{MixedCovThm}
 For general locally finite and diffuse driving measures $\Lambda$ we have
 $$ \Cov({\cal V}_{Y(t,W)},{\cal E}_{Y(t,W)}) = \int_{\bar{[} {\Bbb R}^2 \bar{]}} \left( \frac{1}{2} \Delta^e \otimes \ell(\cdot \cap e) \  {\cal I}^1(s^2 \exp(-s\Lambda([e]));t) \right.$$
 \begin{equation}\label{CovMix5}
  \left.+ \Lambda[\cdot \cap e] \otimes \ell(\cdot \cap e) \ {\cal I}^2(s^2 \exp(-s\Lambda([e]));t) \right)\lparen (\Lambda \times \Lambda)
      \cap \Lambda \rparen(de).
 \end{equation}
\end{theorem}
 Observe, that the first term in the integral (\ref{CovMix5}) takes into account the pairs consisting of a
 vertex constituting the terminus of a maximal edge and the maximal edge itself, whereas the second term
 corresponds to pairs consisting of a vertex lying in the relative interior of a maximal edge and the maximal edge. 
 In analogy to the case of Theorem \ref{CovMeasThm}, vertex-edge pairs where the vertex is
 not adjacent to the edge bring no contribution to the considered covariance structure.

\subsection{Variance Calculation for the Stationary and Isotropic Case}\label{secvariso}

We specialize now the results obtained in the preceding Section \ref{secvargen} to the stationary and isotropic case, i.e. we consider stationary and isotropic random tessellations in the plane that are stable under iteration (STIT tessellations). Up to reparametrization, it means taking the driving measure
$\Lambda_{iso}$ \--- the isometry invariant measure on the space of lines in the plane with length density one.
Recall first that Crofton's formula \cite[Thm. 5.1.1]{SW} for $d=2$ and $k=1$ implies \begin{equation}\Lambda_{iso}([K])={2\over\pi}V_1(K)={1\over\pi}P(K),\label{invmeasure}\end{equation} where $K\subset{\Bbb R}^2$ is a planar convex body with first intrinsic volume $V_1(K)$ and perimeter length $P(K)$. For a line segment $e\subset{\Bbb R}^2$ this is just
\begin{equation}\Lambda_{iso}([e])={2\over \pi}\ell(e),\label{EQLISOELLREL}
\end{equation}
where $\ell(e)$ stands for the length of $e$. It follows from (\ref{invmeasure}) that $$\Sigma_{\Lambda_{iso}[\cdot]}(Y(s,W))={2\over\pi}\Sigma_{\ell}(Y(t,W)).$$ Moreover, in the context of
(\ref{ACONST}) and (\ref{IntersectMeas}) we have $$ \lparen \Lambda_{iso} \cap \Lambda_{iso} \rparen(dx)
= {2\over \pi} dx \ \ \ \ \mbox{ and } \ \ \ \ 
 A_{\Lambda_{iso}[\cdot]}(Y(s,W))={2\over\pi}\Area(W).$$ Recall from Theorem 4 in \cite{ST} that the variance of the total edge length in $W$ of the stationary and isotropic iteration stable random tessellation $Y(u,W)$, $u>0$ simplifies in our particular case to \begin{equation}\nonumber\Var(\Sigma_{\ell}(Y(u,W)))=\pi\int_0^\infty\overline{\gamma}_W(r)\left(1-e^{-{2\over\pi}ur}\right){dr\over r},\label{isovarexpr}\end{equation} where $\overline{\gamma}_W(r)= \int_{{\cal S}_1} \Area(W \cap
 (W + r u)) \nu_1(du)$ is the isotropized set-covariance function of the window $W$,  with $\nu_1$
standing for the uniform distribution on the unit circle ${\cal S}_1$, see \cite{SKM} for the definition of 
$\overline{\gamma}_W(\cdot)$ and Subsection 4.2 in \cite{ST} for further details.
 Combining this with (\ref{ESLA}) and (\ref{VARSI1}), in view of (\ref{EQLISOELLREL}) we are immediately led to the variance formula.
 An alternative method for deriving this formula directly from (\ref{VARS1thm}) in Theorem
\ref{thmvar} in the case $\Lambda=\Lambda_{iso}$ is to use (\ref{invmeasure}) and the important identity for the intersection measure (\ref{SIntersectMeas}), namely
 \begin{equation}\label{ExplicitIntersectMeas}
  \lparen (\Lambda_{iso} \times \Lambda_{iso}) \cap \Lambda_{iso} \rparen (d \overline{xy}) =
  \frac{4dx dy}{\pi^3 ||x-y||}.
 \end{equation}
 The last equation may be established by a twofold application of the affine Blaschke-Petkantschin formula \cite[Thm. 7.2.7]{SW}
 as shown in \cite[Eq. (50)]{ST} in connection with (\ref{EQLISOELLREL}). It follows
\begin{eqnarray} 
\nonumber & & \Var(\Sigma_1(Y(t,W)))\\
\nonumber &=& {t\over\pi}P(W)+{3\over\pi}\Area(W)t^2-2\int_W\int_W{T_3^{\exp}\left(-{2\over\pi}t\left\|x-y\right\|\right)\over\left({2\over\pi}\left\|x-y\right\|\right)^3}{4dxdy\over\pi^3\left\|x-y\right\|}\\
\nonumber &=& {t\over\pi}P(W)+{3\over\pi}\Area(W)t^2-2\pi\int_0^\infty\overline{\gamma}_W(r){T_3^{\exp}\left(-{2\over\pi}tr\right)\over r^3}dr.
\end{eqnarray} 
Computing now $T_3^{\exp}\left(-{2\over\pi}tr\right)$, we arrive at 
\begin{corollary}\label{ISOvarcor} The variance of the number of maximal edges of a stationary
 and isotropic random iteration stable tessellation $Y(t,W)$ is given by
 $$ \Var(\Sigma_1(Y(t,W)))=$$ 
\begin{equation}\label{ISOvarformula}
 {t \over \pi} P(W) +  
 {3\over\pi}\Area(W)t^2+\int_0^\infty\overline{\gamma}_W(r)\left({4t^2\over\pi r}-{4t\over r^2}+{2\pi\over r^3}\left(1-e^{-{2\over\pi}tr}\right)\right)dr.\end{equation}
\end{corollary}
As an example we may consider for $W$ the ball $B_R^2$ in ${\Bbb R}^2$ with radius $R>0$. In this special case, the isotropized set covariance function takes the special form $$\overline{\gamma}_{B_R^2}(r)=2R^2\arccos\left({r\over 2R}\right)-{r\over 2}\sqrt{4R^2-r^2},\ \ \ 0\leq r\leq 2R.$$ Unfortunately, the arising integral cannot further be simplified.\\ \\ We are now interested in the variance asymptotics for a sequence $W_R=R\cdot W$ of growing observation windows, as $R\rightarrow\infty$. To this end, first note that asymptotically, as $R\rightarrow\infty$, we have $$\int_0^{B(R)} \left({4t^2\over\pi r}-{4t\over r^2}+{2\pi\over r^3}\left(1-e^{-{2\over\pi}tr}\right)\right)dr\sim {4\over\pi}t^2\log R $$
as long as $\log B(R) \sim \log R$, where $B(R)$ stands for some upper integration bound depending on $R$. Here and later, we will write $f(R)\sim g(R)$ whenever $\displaystyle\lim_{R\rightarrow\infty}{f(R)\over g(R)}=1$. Now, the relation $$\overline{\gamma}_{W_R}\sim\Area(W_R)=R^2\Area(W),$$ 
valid uniformly for the argument $r = O(R / \log R)$ and $\overline{\gamma}_{W_R} \to 0$ for
$r = \Omega(R \log R)$ (using the standard Landau notation), implies 
\begin{eqnarray}
\nonumber \Var(\Sigma_1(Y(t,W_R))) &\sim & {1\over \pi} t P(W_R) + {3\over\pi}t^2\Area(W_R)+{4\over\pi}t^2\Area(W_R)\log R\\
\nonumber &\sim & {4\over\pi}t^2\Area(W)R^2\log R.
\end{eqnarray}
Summarizing, we have shown
\begin{corollary}\label{ISOasymptvarcor} Asymptotically, as $R\rightarrow\infty$, we have $$\Var(\Sigma_1(Y(t,W_R)))\sim {4\over\pi}\Area(W)t^2R^2\log R$$ and $$\Var(N_v(Y(t,W_R)))\sim {16\over\pi}\Area(W)t^2R^2\log R $$
where $N_v(Y(t,W_R)) \sim 2 \Sigma_1(Y(t,W_R))$ is the number of vertices of $Y(t,W_R).$
\end{corollary} 
\begin{figure}[t]
\begin{center}
 \includegraphics[width=0.45\columnwidth]{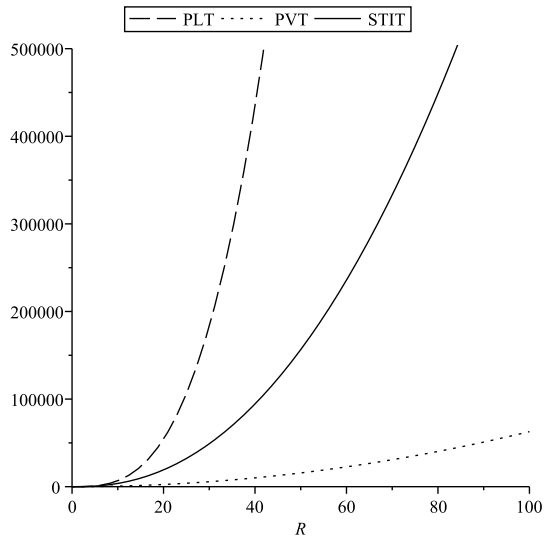}
 \caption{Variance comparison for the number of vertices of Poisson line (PLT), STIT and Poisson-Voronoi (PVT) tessellations for $t=1$ and $W_R=B_R^2$, the ball with radius $R>0$}\label{fig2}
\end{center}
\end{figure}
The formulas show that the geometry of the window $W$ is only reflected by its area in the variance asymptotics. Recall from \cite{HM} that for the Poisson-Voronoi tessellation $\PVT(t,W)$ restricted to some window $W\subset{\Bbb R}^2$ with edge-length density $t>0$ (i.e. the intensity of the underlying Poisson point process equals $t^2\over 4$) we asymptotically have $$\Var(N_v(\PVT(t,W)))\sim 2t^2R^2\Area(W).$$ Moreover, for the stationary and isotropic Poisson line tessellation $\PLT(t,W)$ with length intensity $t$ and restricted to $W$ we have according to \cite{HCLT} $$\Var(N_v(\PLT(t,W)))\sim {4\over\pi^2}t^3R^3\text{CPI}(W;2),$$ where $\text{CPI}(W;2)$ is the $2$-nd-order chord-power integral of $W$ in the sense of \cite[Chap. 8.6]{SW}. The appearance of $\text{CPI}(W;2)$ in the latter formula means that beside the area of $W$, also its shape plays asymptotically an important role. The variances of the different tessellation models are illustrated in Figure \ref{fig2} for a a sequence $W_R=B_R^2$ of circles with radius $R$. In this case it holds $\text{CPI}(B_1^2;2)={16\pi\over 3}$ as can be concluded from Thm. 8.6.6 in \cite{SW} with a corrected constant.\\ In particular, the formulas from the last corollary establish weak long range dependencies (cf. \cite{HS}) for the point process of maximal edge-midpoints and the point process of vertices, since $$\lim_{R\rightarrow\infty}{\Var(\Sigma_1(Y(t,W_R)))\over\Area(W_R)}=\lim_{R\rightarrow\infty}{\Var(N_v(Y(t,W_R)))\over\Area(W_R)}=\infty.$$
As explained at the end of Subsection \ref{secgenPCF}, these dependencies are propagated by long maximal edges on which the vertices are lying and the log-term in the asymptotic variance formula reflects the weakness of these long range dependencies, whereas in contrast to the STIT model, Poisson line tessellations have strong long range dependencies, since $\Area(W_R)^{-1}\Var(N_v(\PLT(t,W_R)))$ grows polynomially (linearly) in $R$, while Poisson-Voronoi tessellations do not have long range dependencies at all. In fact, the maximal edges almost surely have finite length and there are no full straight lines in the tessellation $Y(t)$. In contrast to this, Poisson line tessellations consist by definition only of full lines and the spatial dependencies are in this case much stronger due to the geometric structure of these processes. For Poisson-Voronoi tessellations, we have independence of local geometries whenever the observation regions are far enough from each other, which means that we have 'almost' independence for the point process of vertices, see \cite{HM}. In this sense the STIT tessellations exhibit features intermediate between Poisson-Voronoi and Poisson line tessellations.

\subsection{Vertex Pair-Correlations for the Stationary and Isotropic Case}\label{secPCFiso}

Having made in Section \ref{secgenPCF} very general computations for the covariance measure of the point process of vertices of an iteration infinitely divisible random tessellation $Y(t,W)$, we specialize now -- as in the last Subsection \ref{secvariso} -- to the stationary and isotropic set-up and consider a stationary and isotropic random tessellation $Y(t,W)$ that is iteration stable, i.e. a stationary and isotropic random STIT tessellation with driving measure $\Lambda_{iso}.$ Firstly, recall the relation (\ref{ExplicitIntersectMeas}). Plugging this expression into (\ref{CovMeasExpr}), using in addition that $\Lambda_{iso}([\cdot\cap e])={2\over\pi}\ell(\cdot\cap e)$ and applying the substitution $y=x+u$ yields
\begin{eqnarray}
\nonumber & & \Cov({\cal V}_{Y(t)}) = \int_{{\Bbb R}^2}\int_{{\Bbb R}^2}\left[{1\over 2}
 (\Delta^{\overline{Ou}}\otimes\Delta^{\overline{Ou}}) \circ (\vartheta_x)^{-1} \
 {\cal I}^1(s^2e^{-{2\over\pi}\left\|u\right\|s};t)\right.
\\
\nonumber & & +  \left(\Delta^{\overline{Ou}}\otimes{2\over\pi}\ell(\cdot\cap\overline{Ou})+{2\over\pi}\ell(\cdot\cap\overline{Ou})
\otimes\Delta^{\overline{Ou}}\right) \circ (\vartheta_x)^{-1} \ {\cal I}^2(s^2e^{-{2\over\pi}\left\|u\right\|s};t)\\
\nonumber & & +  
\left.4\left({2\over\pi}\ell(\cdot\cap\overline{Ou})\otimes{2\over\pi}\ell(\cdot\cap\overline{Ou})\right)
\circ (\vartheta_x)^{-1} \ {\cal I}^3(s^2e^{-{2\over\pi}\left\|u\right\|s};t)\right] \ {4dxdu\over\pi^3\left\|u\right\|},
\end{eqnarray}
where we have abbreviated by $O$ the origin and the usual length measure by $\ell$ and where
$\vartheta_x$ stands for the diagonal shift $\vartheta_x(v,w) = (v+x,w+x),\; v,w \in {\Bbb R}^2.$
The covariance measure $\Cov({\cal V}_{Y(t)})$ can be reduced in the sense of \cite[Sec. 8.1]{DVJ}
and, by Proposition 8.1.I(b) there, the reduced covariance measure $\widehat{\Cov}({\cal V}_{Y(t)})$ has the form
\begin{eqnarray}
 \nonumber \widehat{\Cov}({\cal V}_{Y(t)}) &=& \int_{{\Bbb R}^2} \left[{1\over 2} (\delta_u + \delta_{-u} + 2\delta_O)
 \ {\cal I}^1(s^2e^{-{2\over\pi}\left\|u\right\|s};t)
   \right.\\
\nonumber & & +{4\over\pi}\ell(\cdot\cap\overline{(-u)u})  \ {\cal I}^2(s^2e^{-{2\over\pi}\left\|u\right\|s};t)\\
\nonumber & & \left. + 4 \left( {4 \over \pi^2} \int_{\overline{Ou}} \int_{\overline{Ou}} \delta_{v-w} \ell(dv) \ell(dw)
 \right) \ {\cal I}^3(s^2e^{-{2\over\pi}\left\|u\right\|s};t)\right] \ {4du\over\pi^3\left\|u\right\|}
\end{eqnarray}
where we have used the fact that
 $\int \delta_{v-w} (\Delta^{\overline{Ou}} \otimes \Delta^{\overline{Ou}})d(v,w) =
\delta_u + \delta_{-u} + 2 \delta_O$ and $\int \delta_{v-w} (\Delta^{\overline{Ou}} \otimes \ell(\cdot \cap 
 \overline{Ou}))d(v,w) = 2\ell(\cdot \cap \overline{(-u)u}).$
Recall now, see again \cite{DVJ}, that the measure $\widehat{\Cov}({\cal V}_{Y(t)})$ and the reduced second moment measure ${\cal K}({\cal V}_{Y(t)})$ are related by $$\widehat{\Cov}({\cal V}_{Y(t)})={\cal K}({\cal V}_{Y(t)})-\lambda^2\ell_{{\Bbb R}^2},$$ see \cite[Eq. (8.1.6)]{DVJ}, where $\ell_{{\Bbb R}^2}$ is the Lebesgue measure in the plane
and where $\lambda$ stands for the intensity of ${\cal V}_{Y(t)}$. Hence, taking into account that the
vertex intensity $\lambda$ equals ${2\over\pi}t^2$, see \cite{NW06}, and transforming into polar coordinates gives us
\begin{eqnarray}
\nonumber {\cal K}({\cal V}_{Y(t)}) &=& \frac{4}{\pi^3} \int_0^{2\pi} \int_0^{\infty} (\delta_O + \delta_{r e^{i \varphi}})
{\cal I}^{1}(s^2 e^{-{2 \over \pi} r s};t) + {8\over \pi} \ell(\cdot \cap \overline{O re^{i \varphi}})
 {\cal I}^{2}(s^2 e^{-{2 \over \pi} r s};t) \\
 \nonumber & & + 4 \left( \frac{8}{\pi^2} \int_0^r (r-\rho) \delta_{\rho e^{i \varphi}} d\rho\right)  
 {\cal I}^{3}(s^2 e^{-{2 \over \pi} r s};t) dr d\varphi + \left(\frac{2}{\pi} t^2\right)^2 \ell_{\Bbb R}^2.
\end{eqnarray}
From the last expression we can now calculate Ripley's K-function $$K(R) :=\left({\pi\over 2t^2}\right)^2{\cal K}(B_R^2),$$ often also considered in the factorial version $\tilde{K}(R)$ with $$K(R)=\tilde{K}(R)+\left({\pi\over 2t^2}\right)^2{\cal K}(\{0\}),$$ see \cite[Eq. (8.1.12)]{DVJ} or \cite[Chap. 4.5]{SKM}. We obtain 
\begin{eqnarray}
\nonumber K(R) &=& {2\over t^4}\int_0^\infty (1+{\bf 1}[r\leq R])
{\cal I}^1(s^2e^{-{2\over\pi}rs};t)+{8\over\pi}\min(r,R) {\cal I}^2(s^2e^{-{2\over\pi}rs};t)\\
\nonumber & & +
 \frac{32}{\pi^2} \left(r \min(r,R) - {1\over 2}\min(r,R)^2\right)  {\cal I}^3(s^2e^{-{2\over\pi}rs};t)dr+\pi R^2.
\end{eqnarray}
Splitting the integral into two parts, one integral over $[0,R]$ and another over $[R,\infty)$, yields
\begin{eqnarray}
K(R) &=& \pi R^2+{2\over t^4}\int_0^R 2{\cal I}^1(s^2e^{-{2\over\pi}rs};t)+{8\over\pi}r{\cal I}^2(s^2e^{-{2\over\pi}rs};t)\label{kfunctionfinal}\\
 \nonumber & & \hspace{4.6cm}+{16\over\pi^2}r^2{\cal I}^3(s^2e^{-{2\over\pi}rs};t)dr\\
\nonumber &  & +{2\over t^4}\int_R^\infty {\cal I}^1(s^2e^{-{2\over\pi}rs};t)+{8\over\pi}R{\cal I}^2(s^2e^{-{2\over\pi}rs};t)\\
\nonumber & & \hspace{4.6cm} +{32\over\pi^2}\left(rR-{R^2\over 2}\right){\cal I}^3(s^2e^{-{2\over\pi}rs};t)dr.
\end{eqnarray}
Using (\ref{I1}), (\ref{I2}), (\ref{I3}) we finally obtain by using the definition $$g(r)={1\over 2\pi r}{d\over dr}K(r)={1\over 2\pi r}{d\over dr}\tilde{K}(r)$$ of the pair-correlation function, describing the 
normalized vertex density in the distance $r$ from a typical vertex, 
\begin{corollary}\label{CORpcfiso} The pair-correlation function of the vertex point process ${\cal V}_{Y(t)}$ of a stationary and isotropic random STIT tessellation $Y(t)$ with edge-length density $t>0$ equals
\begin{equation}\label{PCFformula}
g(r)=1+{2\over t^2r^2}-{\pi\over t^3r^3}+{\pi^2\over 4t^4r^4}-\left({1\over 2t^2r^2}-{\pi\over 2t^3r^3}+{\pi^2\over 4t^4r^4}\right)e^{-{2\over\pi}tr}.
\end{equation}
\end{corollary}
\begin{figure}[t]
\begin{center}
 \includegraphics[width=0.45\columnwidth]{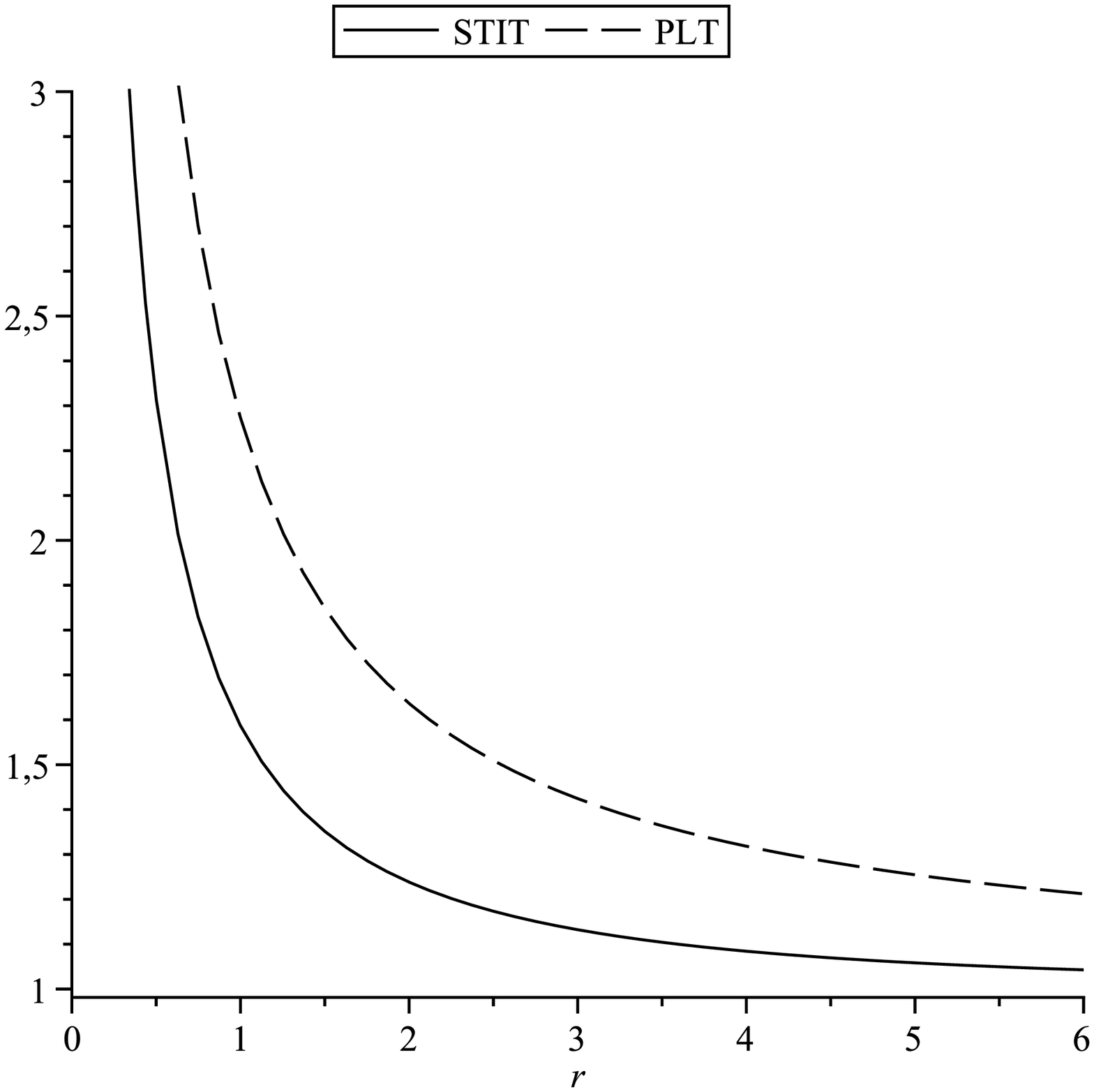}
 \includegraphics[width=0.45\columnwidth]{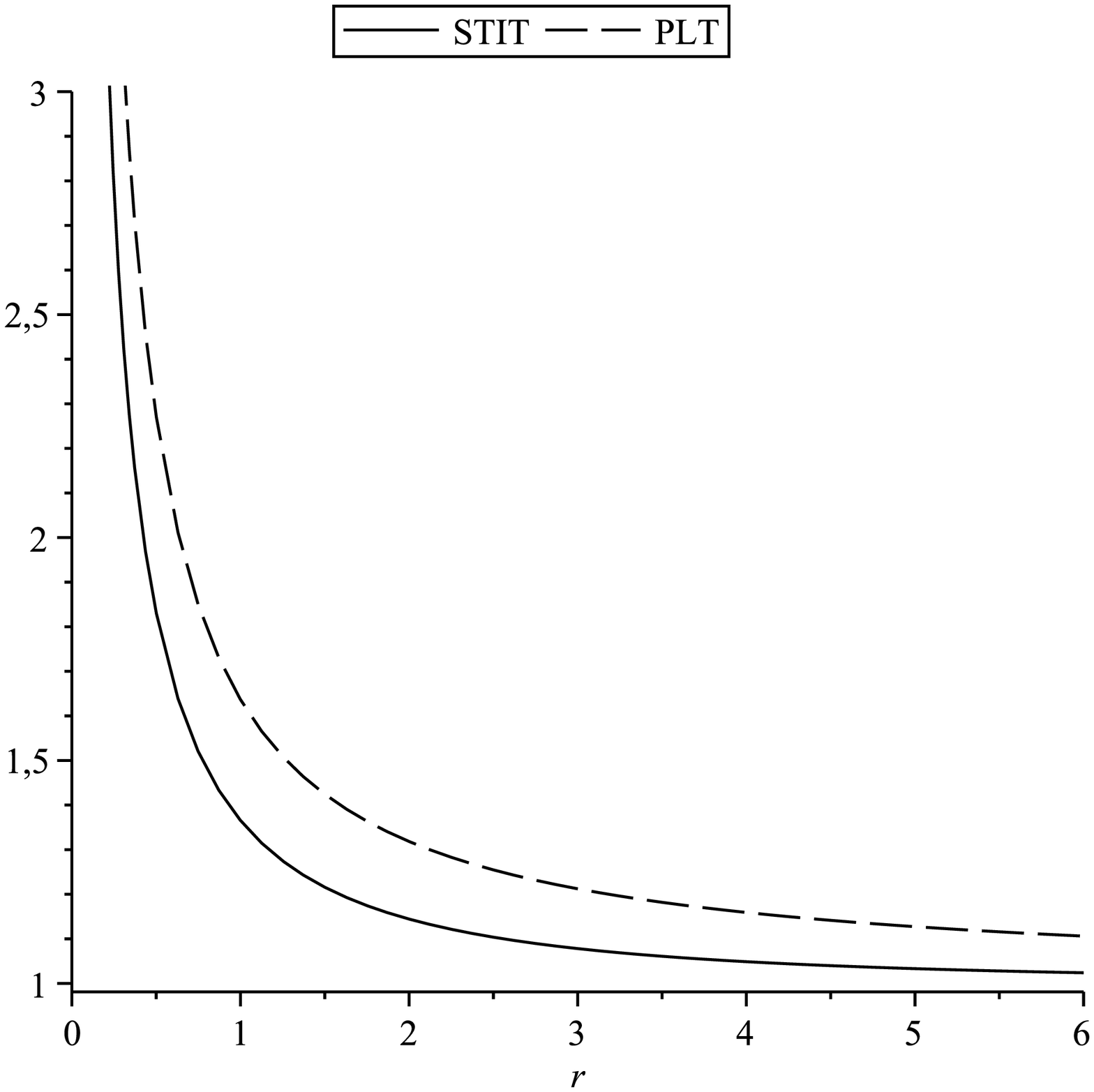}
 \caption{Pair-Correlation function of the point process of vertices of a Poisson line (PLT), STIT and Poisson-Voronoi tessellation (PVT) (left) and the cross-correlation function of the vertex point process and the length measure (right) of a Poisson line (PLT) and a STIT each time with edge-length density $1$}\label{figpcf}
\end{center}
\end{figure}
Note in this context that in the case of a stationary and isotropic Poisson line tessellation with intensity $t>0$, the pair-correlation function of the point process of vertices takes the form $$g^{PLT(t)}(r)=1+{4\over \pi tr},$$ which can easily be concluded from Slivnyak's theorem for Poisson point processes, see \cite{SKM} (here applied to the Poisson process of lines). A comparison of the pair-correlation functions $g(r)$ of the STIT tessellation and $g^{PLT(t)}(r)$ and that of a Poisson-Voronoi tessellation is shown in Figure \ref{figpcf}. However, in contrast to Poisson line and STIT tessellation, the structure of the pair-correlation function of the point process of vertices of a Poisson-Voronoi tessellation is much more complicated. It may be expressed by a sum of integrals of rather involved functions, which cannot be explicitly evaluated. For details and the non-trivial numerical computations we refer to \cite{HM}.\\ It is interesting to compare the pair-correlation formula in Corollary \ref{CORpcfiso} with the
information on the variance of the number of edges provided in Corollary \ref{ISOvarcor}. To this end, use the
variance formula
$$\Var(N_v(Y(t,W)))=2\pi\left({2\over\pi}t^2\right)^2\int_0^\infty\overline{\gamma}_W(r)[g(r)-1]rdr+{2\over\pi}t^2\Area(W), $$
see \cite[Eq. 4.5.7]{SKM} and Girling's formula thereafter, and compare it with the variance expression in Corollary \ref{ISOvarcor}. Taking into account that the number of vertices is, modulo boundary effects, twice the number of edges, we
should have agreement of the leading term $\frac{4t^2}{\pi r}$ in the integral in (\ref{ISOvarformula})
with ${1\over 4}\cdot 2\pi r\left({2\over\pi}t^2\right)^2$ times the leading term $\frac{2}{t^2 r^2}$ of
$g(r)-1$ in (\ref{PCFformula}), where ${1\over 4}$ comes from switching between edge and vertex
counts, the factor $2\pi r$ comes from transformation into polar coordinates and the remaining factor is 
the squared intensity of the vertex point process. Comparing these expressions we readily obtain
the required agreement of leading terms determining the prefactor in the $O(R^2 \log R)$-variance
asymptotics. The remaining lower order $o(1/r^2)$ terms in the pair correlation function (\ref{PCFformula})
do not have to and do not agree with their suitably normalized counterparts in (\ref{ISOvarformula}), because
the latter contains additional area order corrections and, moreover, takes into account the boundary effects
caused by edges hitting the boundary $\partial W_R$ of $W_R$.\\
Another aspect that can be compared concerns the radial distribution function. For a stationary and isotropic random point process in the plane with intensity $\lambda>0$ and K-function $K(r)$ the radial distribution function $\rho(r)$ is defined by $$\rho(r)=\lambda{dK(r)\over dr}.$$ Writing from now on $\rho(r)$ for the radial distribution function of the point process ${\cal V}_{Y(t)}$ of vertices of a stationary and isotropic STIT tessellation $Y(t)$ with edge-length density $t>0$ we can use (\ref{kfunctionfinal}) to conclude $$\rho(r)=4t^2r+{8\over r}-{4\pi\over tr^2}+{\pi^2\over t^2r^3}-\left({2\over r}-{2\pi\over tr^2}+{\pi^2\over t^2r^3}\right)e^{-{2\over\pi}tr},$$ whereas for the radial distribution function $\rho^{\PLT(t)}(r)$ of the vertex point process of a stationary and isotropic Poisson line tessellation with edge-length density $t>0$ we obtain $$\rho^{\PLT(t)}(r)=2t^2r+{8t\over\pi}$$ from Slivnyak's theorem. It means that asymptotically we have $$\rho(r)\sim 2\rho^{\PLT(t)}(r),\ \ \ \text{as}\ \ \ r\rightarrow\infty.$$
\begin{remark} There are different normalizations available for the reduced second-moment measure $\cal K$ and Ripley's K-function in the existing literature as for example \cite{DVJ} or \cite{SKM}. We decided here not to normalize $\cal K$ by one over the squared intensity, $1\over\lambda^2,$ but we normalize the K-function by that factor in order ensure that the pair-correlation function $g(r)$ tends to $1$ as $r\rightarrow\infty$. This is done to keep the formulas consistent with those from previous papers on STIT tessellations. This convention will also be adopted in the next subsection.
\end{remark}

\subsection{Edge-Vertex Correlations in the Stationary and Isotropic Case}\label{seccrosspcfiso}

Our interest here is focused on the cross-covariance measure $\Cov({\cal V}_{Y(t)},{\cal E}_{Y(t)})$ of a stationary and isotropic random STIT tessellation $Y(t)$ with edge-length density $t>0$ and driving measure $\Lambda_{iso}.$ It describes the correlations between the stationary and isotropic random point process of vertices and the stationary and isotropic random length measure concentrated on the edges of $Y(t)$. The study of this measure was proposed in \cite{SO1} and \cite{SO2} and we recall some general definitions now. Let $\Phi_1$ and $\Phi_2$ be stationary and isotropic random measures in ${\Bbb R}^2$ with respective intensities $\lambda_1>0$ and $\lambda_2>0$. For a Borel set $B\subset{\Bbb R}^2$ we introduce the random measures $${\cal K}_{12}(B):={\Bbb E}\int_{[0,1]^2}\Phi_2(B+x)\Phi_1(dx),\ \ \ {\cal K}_{21}(B):={\Bbb E}\int_{[0,1]^2}\Phi_1(B+x)\Phi_2(dx).$$ The measure ${\cal K}_{12}$ describes $\Phi_2$ as seen from the typical point of $\Phi_1$ and ${\cal K}_{21}$ describes the measure $\Phi_1$  regarded from the typical point of $\Phi_2$ in the sense of Palm distributions. In \cite{SO2} it was shown that ${\cal K}_{12}(B)={\cal K}_{21}(-B),$ which in particular implies for the ball $B_r^2$ with radius $r>0$ the identity ${\cal K}_{12}(B_r^2)={\cal K}_{21}(B_r^2)$. The cross-K-function $K_{12}(r)=K_{21}(r)$ may now be introduced as $$K_{12}(r):={1\over\lambda_1\lambda_2}{\cal K}_{12}(B_r^2)$$ and the cross-correlation function $g_{12}(r)$ of the random measures $\Phi_1$ and $\Phi_2$ is defined by \begin{equation}g_{12}(r):={1\over 2\pi r}{dK_{12}(r)\over dr},\label{cpcfgen}\end{equation} compare with \cite{SO1} an \cite{SO2}. Informally, we could say that $\lambda_1K_{21}(r)$ or $\lambda_2K_{12}(r)$ is the expectation of $\Phi_1(B_r^2+x)$ or $\Phi_2(B_r^2+x)$ at the typical point $x$ of $\Phi_2$ or $\Phi_1$, respectively.\\ The theory is now applied to our setting and we take for $\Phi_1$ the point process ${\cal V}_{Y(t)}$ and for $\Phi_2$ the random measure ${\cal E}_{Y(t)}$ (recall the definitions from Section \ref{secgenEVPCF}). We have in this special situation $\lambda_1 = \frac{2}{\pi} t^2$
and $\lambda_2 = t,$ see \cite{NW06}.
We can now use Theorem \ref{MixedCovThm} together with (\ref{ExplicitIntersectMeas}) and 
(\ref{invmeasure}) to obtain, under the substitution $y = x + u$,
\begin{eqnarray}
\nonumber \Cov({\cal V}_{Y(t)},{\cal E}_{Y(t)}) &=& \int_{{\Bbb R}^2} \int_{{\Bbb R}^2}\left[ \frac{1}{2}
 \Delta^{\overline{Ou}} \otimes \ell(\cdot \cap \overline{O u}) {\cal I}^1(s^2 e^{-\frac{2}{\pi} ||u||s};t) \right.\\
\nonumber & & \left.+\frac{2}{\pi} \ell(\cdot \cap \overline{O u})^{\otimes 2} \ 
        {\cal I}^2(s^2 e^{-\frac{2}{\pi} ||u||s};t) \right] \circ \vartheta_x^{-1}\frac{4du dx}{\pi^3 ||u||}.
\end{eqnarray}
Consequently, we end up with
\begin{eqnarray}
\nonumber {\cal K}_{12} &=& \int_{{\Bbb R}^2} \left[ \ell(\cdot \cap \overline{O u}) {\cal I}^1(s^2 e^{-\frac{2}{\pi} ||u||s};t)\right.\\
\nonumber & & \left.+ \frac{2}{\pi} \int_{\overline{Ou}} \int_{\overline{Ou}} \delta_{v - w} \ell(dv) \ell(dw) 
        {\cal I}^2(s^2 e^{-\frac{2}{\pi} ||u||s};t) \right]\frac{4du dx}{\pi^3 ||u||}+{2\over\pi}t^3\ell_{{\Bbb R}^2},
\end{eqnarray}        
whereby, upon passing to polar coordinates, we have
\begin{eqnarray}
\nonumber {\cal K}_{12} &=& \frac{4}{\pi^3} \int_0^{\infty} \int_0^{2\pi} \ell(\cdot \cap \overline{O r e^{i\varphi}}) {\cal I}^1(s^2 e^{-\frac{2}{\pi} rs};t)\\
\nonumber & & + \frac{4}{\pi} \left( \int_0^r (r-\rho) \delta_{\rho e^{i \varphi}} d\rho \right) 
        {\cal I}^2(s^2 e^{-\frac{2}{\pi} rs};t) d\varphi dr +{2\over\pi}t^3\ell_{{\Bbb R}^2}.
\end{eqnarray}
Recalling the definition of $K_{12}$ and using again
that in our set-up $\lambda_1 = \frac{2}{\pi} t^2$ and $\lambda_2 = t,$ 
we obtain
\begin{eqnarray}
\nonumber K_{12}(R) &=& \pi R^2+\frac{4}{\pi t^3} \int_0^{\infty} \min(r,R) \ {\cal I}^1(s^2 e^{-\frac{2}{\pi} rs};t)\\
\nonumber & & + \frac{4}{\pi} \left(r\min(r,R) - {1\over 2}\min(r,R)^2\right) \ {\cal I}^2(s^2 e^{-\frac{2}{\pi} rs};t) dr\\
\nonumber &=& \pi R^2+{4\over\pi t^3}\int_0^Rr{\cal I}^1(s^2 e^{-\frac{2}{\pi} rs};t)+{2\over\pi}r^2{\cal I}^2(s^2 e^{-\frac{2}{\pi} rs};t)dr\\
\nonumber & & +{4\over\pi t^3}\int_R^\infty R{\cal I}^1(s^2 e^{-\frac{2}{\pi} rs};t)+{4\over\pi}\left(rR-{R^2\over 2}\right){\cal I}^2(s^2 e^{-\frac{2}{\pi} rs};t)dr.
\end{eqnarray}
Using now (\ref{cpcfgen}) together with (\ref{I1}), (\ref{I2}) and (\ref{I3}) we arrive at
\begin{corollary} The cross-correlation function of the vertex point process ${\cal V}_{Y(t)}$ and the random length measure ${\cal E}_{Y(t)}$ of a stationary and isotropic random STIT tessellation $Y(t)$ with edge-length density $t>0$ equals $$g_{12}(r)=1+{1\over t^2r^2}-{\pi\over 4t^3r^3}-\left({1\over 2t^2r^2}-{\pi\over 4t^3r^3}\right)e^{-{2\over\pi}tr}.$$
\end{corollary}
In contrast to this formula, the same cross-correlation function $g_{12}^{\PLT(t)}(r)$ for a stationary and isotropic Poisson line tessellation with edge-length density $t>0$ is given by $$g_{12}^{\PLT(t)}(r)=1+{2\over\pi tr},$$ which can easily be obtained from Slivnyak's theorem for which we refer to \cite{SKM}. A comparison of both functions is shown in Figure \ref{figpcf}. We would like to point out that the corresponding cross-correlation function is unknown until know for the Poisson-Voronoi model.

\section{Central Limit Theory}\label{secclt}

 In this section we will study the functional central limit problem for the total edge count and edge length
 processes induced by a STIT tessellation in growing windows $W_R = R W,\; R \to \infty,$ with $W$
 standing for some compact convex set of non-empty interior, to remain fixed throughout the section.  
 We assume that the measure $\Lambda$ is translation invariant, i.e. we are in the STIT regime.\\
 Define the rescaled total edge length process
\begin{equation}\nonumber\label{LGTHPROC}
 {\cal L}^{R,W}_t := \frac{1}{R \sqrt{\log R}} \bar\Sigma_{\Lambda_{iso}[\cdot]}(Y(t+1/\log R,W_R)),
\end{equation}
 as well as the rescaled total edge count process
\begin{equation}\nonumber\label{ECTPROC}
 {\cal C}^{R,W}_t := \frac{1}{R \sqrt{\log R}} \bar\Sigma_1(Y(t+1/\log R,W_R)),\; t \in [0,1].
\end{equation}
 The main result of this section is
 \begin{theorem}\label{FuncCLT}
  The processes $({\cal L}^{R,W}_s,{\cal C}^{R,W}_s)_{s \in [0,1]}$ converges jointly in law, as $R\to\infty,$
  on the space ${\cal D}([0,1];{\Bbb R}^2)$ of ${\Bbb R}^2$-valued c\`adl\`ag functions on $[0,1]$
  endowed with the usual Skorokhod $J_1$-topology, \cite[Chap. VI.1]{JS} or 
  \cite[Chap. 3, Sec. 14]{BIL}, to the process $t \mapsto (\xi,t\xi),$ where
  $\xi$ is a normal random variable with variance $V(\Lambda[\cdot],W)$ which is given by (63)
  in \cite{ST} for $\phi = \Lambda[\cdot]$ there.
\end{theorem}
 Rather than giving a general formula for $V(\phi,W)$, we refer the reader to Proposition 1 and (72) of \cite{ST},
 where the general case is considered and $V(\phi,W)$ is expressed as a weighted mean width of an associated zonoid and its polar body. Here, we only mention the fact that for the particular
 isotropic case we simply have
 $$ V(\Lambda_{iso}[\cdot],W) = {4\over \pi}\Area(W), $$ 
 see the discussion following Proposition 1 in \cite{ST}.\\ 
 The phenomenon observed in Theorem \ref{FuncCLT} above deserves a short discussion. Namely,
 although both $\bar\Sigma_{\Lambda[\cdot]}(Y(t,W_R))$ and $\bar\Sigma_1(Y(t,W_R))$ exhibit
 fluctuations of the order $R \log R,$ the mechanisms in which these fluctuations arise are of
 a rather different nature:
 \begin{itemize}
  \item  As shown in \cite[Thm. 6]{ST}, 
           the leading-order deviations of $\bar\Sigma_{\Lambda[\cdot]}$ arise very early in the course
           of the MNW-construction, in its initial stages usually referred to as the {\it big bang} phase.
           Here, this is the time period $[0,1/\log R].$ During the later stages of the construction,
           i.e. the time interval $(1/\log R,1],$ the variance increase is of lower order and any newly
           arising fluctuations are negligible compared to those originating from the {\it big bang}.
           In the asymptotic picture this means the initial fluctuation remains {\it frozen} throughout the
           rest of the dynamics, whence the constant limit for ${\cal L}^{R,W}_t$ (note at this point
           that the Brownian limit for the length process obtained in Theorem 6 in \cite{ST} referred
           to a different time flow).
 \item In contrast, the deviations of $\bar\Sigma_1$ arise and cumulate constantly in time $t$ with rate proportional
          to $t$ times the initial {\it big bang} fluctuation of  $\bar\Sigma_{\Lambda[\cdot]}.$ Thus, as opposed
          to that of $\bar\Sigma_{\Lambda[\cdot]},$ the variance of $\bar\Sigma_1$ exhibits a non-vanishing
          quadratic dependency on $t$ even in large $R$ asymptotics.
 \end{itemize}
 Thus, in large $R$ asymptotics, we have the following intuitive picture: denoting by $\Xi$ the initial
 {\it big bang} fluctuation of $\Sigma_{\Lambda[\cdot]}$  we can effectively use the following
 first-order approximations: $\bar\Sigma_{\Lambda[\cdot]}(Y(t,W_R)) \approx \Xi$ and
 $\bar\Sigma_1(Y(t,W_R)) \approx t \Xi$ valid for $t \in (1/\log R,1].$ 

\paragraph*{Proof of Theorem \ref{FuncCLT}}
 Consider the auxiliary process
 \begin{equation}\nonumber\label{SIGMAHAT}
  \hat\Sigma_1(Y(t,W)) := \bar\Sigma_1(Y(t,W)) - \int_0^t \bar{A}_1(Y(s,W)) ds,
 \end{equation}
 which is in view of (\ref{A1}) the same as
 $$ \bar\Sigma_1(Y(t,W)) - \int_0^t \bar\Sigma_{\Lambda[\cdot]}(Y(s,W)) ds $$
 and which is a centered $\Im_t$-martingale by (\ref{ORD1mart}).  Squaring and
 taking expectations we get
 \begin{eqnarray}
 \nonumber {\Bbb E}(\hat\Sigma_1(Y(t,W)))^2 &=& \Var(\Sigma_1(Y(t,W))) - 2 
     \int_0^t {\Bbb E}[\bar\Sigma_1(Y(t,W)) \bar\Sigma_{\Lambda[\cdot]}(Y(s,W))] ds\\
 & & +\ {\Bbb E} \int_0^t \int_0^t \bar\Sigma_{\Lambda[\cdot]}(Y(s,W)) 
         \bar\Sigma_{\Lambda[\cdot]}(Y(u,W)) du ds.
\label{VARSIGHAT1}
 \end{eqnarray}
 Using that for $s < t$
 $$ {\Bbb E}(\bar\Sigma_{\Lambda[\cdot]}(Y(t,W))|\Im_s) = \bar\Sigma_{\Lambda[\cdot]}(Y(s,W)) $$
 and
 $$ {\Bbb E}(\bar\Sigma_1(Y(t,W))|\Im_s) - \bar\Sigma_1(Y(s,W)) = $$ $$ \int_s^t 
 {\Bbb E}(\bar\Sigma_{\Lambda[\cdot]}(Y(u,W))|\Im_s) du =
 (t-s) \bar\Sigma_{\Lambda[\cdot]}(Y(s,W)), $$
 as follows by the martingale property of $\bar\Sigma_{\Lambda[\cdot]}(Y(t,W))$
 and $\hat{\Sigma}_1(Y(t,W)),$  we can rewrite (\ref{VARSIGHAT1}) by taking first conditional expectations as
 \begin{eqnarray}
 \nonumber & & {\Bbb E}(\hat\Sigma_1(Y(t,W)))^2\\
 \nonumber &=& \Var(\bar\Sigma_1(Y(t,W))) - 2 
     \int_0^t (t-s) {\Bbb E}\bar\Sigma^2_{\Lambda[\cdot]}(Y(s,W)) ds\\
 \nonumber &-& 2 \int_0^t  {\Bbb E}\bar\Sigma_1(Y(s,W)) \bar\Sigma_{\Lambda[\cdot]}(Y(s,W)) ds +
       2  \int_0^t \int_0^s {\Bbb E} \bar\Sigma^2_{\Lambda[\cdot]}(Y(u,W))
     du ds
 \end{eqnarray}
 and hence, with the second and fourth term in the right hand side canceling out, $${\Bbb E}(\hat\Sigma_1(Y(t,W)))^2 =$$
 \begin{equation}\label{VarReduced}
   \Var(\Sigma_1(Y(t,W))) -
  2 \int_0^t  \Cov(\Sigma_1(Y(s,W)) \Sigma_{\Lambda[\cdot]}(Y(s,W))) ds.
 \end{equation}
 The relation (\ref{VarReduced}) combined with (\ref{VARFOR1}) and (\ref{ESLA}) readily yields
 \begin{equation}\label{VarReduced2}
  {\Bbb E}(\hat\Sigma_1(Y(t,W)))^2 = t \Lambda([W]) +  \frac{t^2}{2} \lparen \Lambda \cap \Lambda \rparen(W).
 \end{equation}
 To proceed, define now the auxiliary process
 \begin{equation}\label{PROCDEF2}
  \hat{\cal C}^{R,W}_t := \frac{1}{R \sqrt{\log R}} \hat{\Sigma}_1(Y(t+1/\log R,W_R)) = 
  {\cal C}^{R,W}_t  - \int_{-1/\log R}^t  {\cal L}^{R,W}_s ds.
 \end{equation}
  Using that $\Lambda([W_R]) = O(R)$ and $\lparen \Lambda \cap \Lambda \rparen(W_R) = O(R^2)$
  we conclude from (\ref{VarReduced2}) that
  $$ \lim_{R\to 0} {\Bbb E}(\hat{\cal C}^{R,W}_1)^2 = 0. $$
  Thus, Doob's $L^2$-maximal inequality \cite[Thm 3.8(iv)]{KS} implies
  \begin{equation}\label{L2conv}
   \lim_{R\to\infty} {\Bbb E} \sup_{t \in [0,1]} (\hat{\cal C}^{R,W}_t)^2 \leq  \lim_{R\rightarrow\infty}4{\Bbb E}(\hat{{\cal C}}_1^{R,W})^2= 0.
  \end{equation}
 We are now in a position to apply Theorem 6 in \cite{ST} (taking into account that the
 notation ${\cal L}^{R,W}$ used in \cite{ST} corresponds to a time change of ${\cal L}^{R,W}$ 
 as defined here)  to conclude that, as $R \to \infty,$ 
 \begin{equation}\nonumber\label{ConvGauss}
  ({\cal L}^{R,W}_t)_{t \in [0,1]} \Longrightarrow (\xi)_{t \in [0,1]},
 \end{equation}
 that is to say the process $({\cal L}^{R,W}_t)_{t \in [0,1]}$
 converges in law in the Skorokhod space ${\cal D}([0,1];{\Bbb R})$ to the constant process $t \mapsto \xi,$ where, recall,
 $\xi$ is a centered normal random variable with variance $V(\Lambda[\cdot],W).$ Using (\ref{L2conv})
 and recalling the definition (\ref{PROCDEF2}) we see now that the processes 
 $({\cal L}^{R,W}_t,{\cal C}^{R,W}_t)_{t \in [0,1]}$ converge jointly in law
 in ${\cal D}([0,1];{\Bbb R}^2)$ to the process $$t \mapsto\left( \xi,\int_0^t \xi ds\right) = (\xi,t\xi),$$
 which completes the proof of Theorem \ref{FuncCLT}.\hfill $\Box$
 
 \begin{remark} For reasons discussed in detail in Remark 6 in
 \cite{ST}, we expect the rate of the convergence in Theorem \ref{FuncCLT} to be rather slow.
 \end{remark}

\section*{Acknowledgements}
The authors would like to thank Joachim Ohser (Darmstadt) for providing the two pictures of simulated STIT tessellations. We are also grateful to the anynomous reviewers, whose remarks were very helpful in improving the style of this work.\\ The first author was supported by the Polish Minister of Science and Higher Education grant N N201 385234 (2008-2010) and the second author was supported by the Swiss National Science Foundation grant PP002-114715/1.


\begin{thebibliography}{30}\small

\bibitem{BIL} 
 {\sc Billingsley, P.} (1968).
 {\textit{Convergence of Probability Measures.}} Wiley.

\bibitem{DVJ}
 {\sc Daley, D.J.; Vere-Jones, D.} (2003).
 {\textit{An Introduction to the Theory of Point Processes.}} Volume I, Springer.

\bibitem{HCLT}
{\sc Heinrich, L.} (2009). {Central limit theorems for motion-invariant Poisson hyperplanes in expanding convex windows.} \textit{Rendiconti del circolo matematico di Palermo}, Series II, Suppl. \textbf{81} , 187--212.

\bibitem{HM}
{\sc Heinrich, L.; Muche, L.} (2008). {Second-order properties of the point process of nodes in a stationary Voronoi tessellation.} \textit{Math. Nachr.} \textbf{281}, 350-375.

\bibitem{JS}
 {\sc Jacod, J., Shiryaev, A.N.} (2003).
 {\textit{Limit Theorems for Stochastic Processes.}} Grundlehren der Mathematischen Wissenschaften {\bf 288}, Second Ed., Springer. 

\bibitem{KS}
 {\sc Karatzas, I., Shreve, S.E.} (1998). {\textit{Brownian Motion and Stochastic Calculus.}} Graduate Texts in Mathematics,
 Second Ed., Springer.

\bibitem{MNW}
{\sc Mecke, J.; Nagel, W.; Wei\ss , V.} (2008). {A global construction of homogeneous random planar tessellations that are stable under iteration.} \textit{Stochastics} \textbf{80}, 51--67.

\bibitem{NW05}
{\sc Nagel, W.; Wei\ss , V.} (2005). {Crack STIT tessellations: characterization of stationary random tessellations stable with respect to iteration.} \textit{Adv. Appl. Probab.} \textbf{37}, 859--883.

\bibitem{NW06}
{\sc Nagel, W.; Wei\ss , V.} (2006). {STIT tessellations in the plane.} \textit{Rendiconti del circulo matematico di Palermo}, Serie II, Suppl. \textbf{77}, 441--458.

\bibitem{HS}
{\sc Schmidt, H.} (2008). {\textit{Asymptotic Analysis of Stationary Random Tessellations.}} VDM Verlag Dr. M\"uller, Saarbr\"ucken.

\bibitem{SW}
{\sc Schneider, R.; Weil, W.} (2008). {\textit{Stochastic and Integral Geometry.}} Springer.

\bibitem{ST}
{\sc Schreiber, T.; Th\"ale, C.} (2010). {Typical geometry, second-order properties and central limit theory for iteration stable tessellations.} \textit{arXiv: 1001.0990 [math.PR]}.

\bibitem{SKM}
{\sc Stoyan, D.; Kendall, W.S.; Mecke, J.} (1995). {\textit{Stochastic Geometry and its Applications.}} Second Ed., Wiley.

\bibitem{SO1}
{\sc Stoyan, D.; Ohser, J.} (1982). {Correlations between planar random structures with an ecological application.} \textit{Biom. J.} \textbf{24}, 631--647.

\bibitem{SO2}
{\sc Stoyan, D.; Ohser, J.} (1985). {Cross-correlation measures for weighted random measures and their estimation.} \textit{Theory Probab. Appl.} \textbf{29}, 345-355.

\bibitem{NOW}
{\sc Wei\ss , V.; Ohser, J.; Nagel, W.} (2010). {Second moment measure and $K$-function for planar STIT tessellations.}
\textit{Image Anal. Stereol.} \textbf{29}, 121--131.

\end{thebibliography}
\end{document}